\documentclass[11pt,reqno]{amsart}
\usepackage{amsmath,amssymb,amsfonts,amsthm,epsfig,picinpar,xspace}
\usepackage[round]{natbib}
\newcommand{\Z}{{\mathbb Z}}

\newtheorem{theorem}{Theorem}
\newtheorem{lemma}[theorem]{Lemma}
\newtheorem{proposition}[theorem]{Proposition}

\newtheorem{definition}{Definition}

\renewcommand{\Z}{{\mathbb Z}}
\newcommand{\old}[1]{}
\newcommand{\gap}{\gamma}
\newcommand{\hl}{\ensuremath{\check{h}}\xspace}
\newcommand{\hu}{\ensuremath{\hat{h}}\xspace}
\newcommand{\lrs}{Luby, Randall, and Sinclair\xspace}
\newcommand{\dsc}{Diaconis and Saloff-Coste\xspace}

\newcommand{\bd}{Bubley and Dyer\xspace}

\newcommand{\kdk}{Karzanov-Khachiyan\xspace}

\newcommand{\mc}{Markov chain\xspace}

\newcommand{\sref}[1]{\S~\ref{sec:#1}}
\newcommand{\tref}[1]{Theorem~\ref{thm:#1}}
\newcommand{\lref}[1]{Lemma~\ref{lem:#1}}
\newcommand{\fref}[1]{Figure~\ref{fig:#1}}
\newcommand{\eref}[1]{Equation~\eqref{eqn:#1}}
\newcommand{\tbref}[1]{Table~\ref{tbl:#1}}
\newcommand{\U}{{\sf U}\xspace}
\newcommand{\D}{{\sf D}\xspace}
\newcommand{\F}{{\sf F}\xspace}
\newcommand{\p}[1]{{\sf #1}}

\newcommand{\pp}[2]{$\genfrac{}{}{-5pt}{1}{\sf #1}{\sf #2}$}

\newcommand{\df}{\Delta\Phi}
\newcommand{\fm}{\Phi_{\min}}
\newcommand{\fx}{\Phi_{\max}}
\newcommand{\ft}{\Phi_t}
\newcommand{\fp}{\Phi_{t+1}}

\newcommand{\Var}{\operatorname{Var}}
\newcommand{\Cov}{\operatorname{Cov}}
\newcommand{\wgap}{\operatorname{wgap}}
\newcommand{\erf}{\operatorname{erf}}
\newcommand{\TV}{\operatorname{TV}}
\newcommand{\1}{{\hat{1}}}
\newcommand{\0}{{\hat{0}}}
\newcommand{\remark}{\noindent\textit{Remark:}\hspace*{1ex}}
\newcommand{\nib}{\noindent$\bullet$\ \ }
\newcommand{\mapprox}{\stackrel{?}{\approx}}
\newcommand{\meq}{\stackrel{?}{=}}

\begin{document}
\title[Tiling and shuffling mixing times]{\vspace*{-50pt}
Mixing times of lozenge tiling and\\
card shuffling Markov chains}
\author[David B. Wilson]{David Bruce Wilson}
\address{Microsoft Research\\One Microsoft Way\\Redmond, WA 98052\\U.S.A.}
\email{dbwilson@microsoft.com}
\urladdr{http://dbwilson.com}
\thanks{The research that led to this article was done while the author was at MIT, DIMACS at Rutgers, Stanford, and Microsoft Research, and was supported in part by an NSF-postdoctoral fellowship.}
\date{}

\vspace*{-12pt}
\begin{abstract}
We show how to combine Fourier analysis with coupling arguments to bound
the mixing times of a variety of Markov chains.  The mixing time is the number
of steps a Markov chain takes to approach its equilibrium
distribution.  One application is to a class of Markov chains
introduced by \lrs to generate random tilings of regions by lozenges.
For an $\ell\times\ell$ region we bound the mixing time by $O(\ell^4\log
\ell)$, which improves on the previous bound of $O(\ell^7)$, and we show
the new bound to be essentially tight.  In another application we
resolve a few questions raised by \dsc by lower bounding the mixing
time of various card-shuffling Markov chains.  Our lower bounds are
within a constant factor of their upper bounds.  When we use our
methods to modify a path-coupling analysis of \bd, we obtain an $O(n^3\log n)$
upper bound on the mixing time of the \kdk Markov chain for linear
extensions.
\end{abstract}

\maketitle

\section{Introduction}

Using a simple idea, we obtain improved upper and lower bounds on the
mixing times of a number of previously studied Markov chains:

\begin{figwindow}[0,r,%
{\psfig{figure=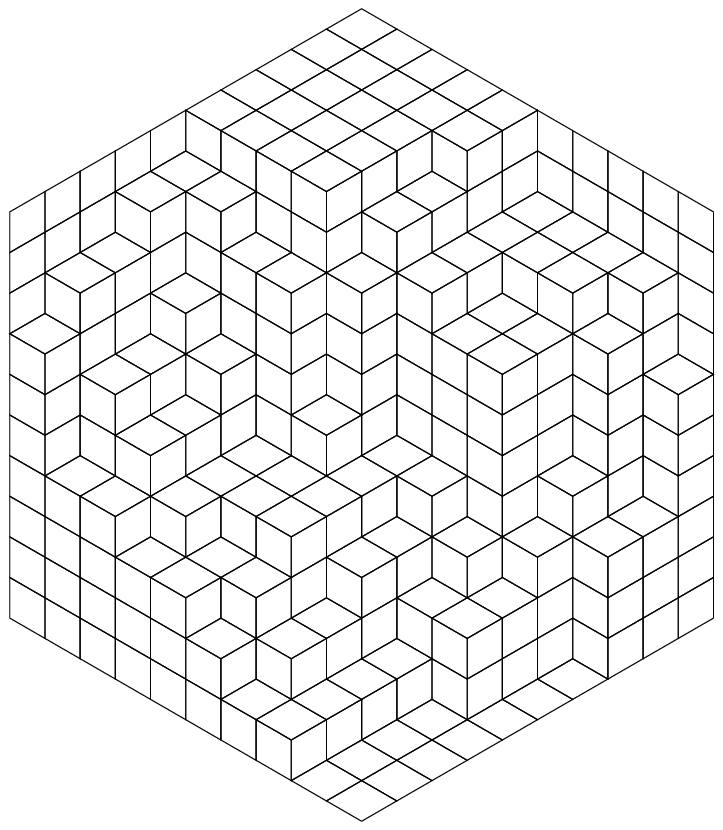}},%
{Random lozenge tiling of the order 10 hexagon, chosen uniformly at random from the 9265037718181937012241727284450000 possible such tilings. Here $\ell=10$, $w=20$, and $n=300$.\label{fig:cube}}]
\nib
A lozenge is a rhombus with angles $120^\circ$ and $60^\circ$ and
sides of unit length.  \fref{cube} shows a random lozenge tiling of a
hexagon.  Random lozenge tilings were originally studied in physics as
a model for dimer systems, and have recently served as an exploratory
tool by people in combinatorics.  \sref{lozenge} gives further back-
ground.  \citet*{luby-randall-sinclair:markov-lattice} proposed
a Markov chain to generate random lozenge tilings of regions,
proved that it runs in time $O(n^4)$ when there are $n$
lozenges, and in later unpublished work reduced the bound to $O(n^{3.5})$.
(They also analyzed domino-
tiling and Eulerian-orientation Markov chains.)
For the regular hexagon with side length $\ell$, for example,
their methods give a bound of
$O(\ell^7)$.  We show here that the
correct mixing time of this Markov chain is $\Theta(\ell^4 \log
\ell)=\Theta(n^2\log n)$, by showing the state to be very far from stationarity
after $\sim(8/\pi^2) \ell^4\log \ell$ steps, and very close to stationarity after
$\sim(48/\pi^2) \ell^4\log \ell$ steps.  The correct constant appears to
be $16/\pi^2$.  For general regions of size $n$ and
width $w$, our upper bound is $\sim(3/\pi^2)w^2 n\log n$.

\nib
Consider the following shuffle on a deck of $n$ cards: with
probability $1/2$ do nothing, otherwise transpose a random adjacent
pair of cards.  How many of these operations are needed before the
deck becomes well-shuffled?  In the years since
\citet*[sect.~4]{aldous:group-walks} showed that $O(n^3 \log n)$
shuffles are enough and that $\Omega(n^3)$ shuffles are necessary to
randomize the deck, there have been a couple of heuristic arguments
\citep{aldous:personal} \citep{diaconis:personal} for why $\Omega(n^3
\log n)$ shuffles should be necessary, but unruly technical
difficulties prevented a rigorous proof from being written down.
Using our method these technical difficulties vanish.  With little
more than algebra and trigonometry, we show that $(1/\pi^2-o(1)) n^3
\log n$ shuffles are necessary to begin to randomize the deck, and
that $(2/\pi^2+o(1)) n^3\log n$ shuffles are enough.  The best
previous published explicit upper bound was $(4+o(1)) n^3 \log n$ shuffles.
The correct constant appears to be $1/\pi^2$.
\end{figwindow}

\nib
We lower bound the mixing time of a few of other shuffles analyzed by
\citet*{diaconis-saloff-coste:comparison}.  They considered
shuffles for which the cards appear on the vertices of a graph, and
the shuffle picks a random edge and transposes the cards at its
endpoints.  They obtained upper bounds on the mixing time for shuffles
based on the $\sqrt{n}\times\sqrt{n}$ grid and on the $(\log_2
n)$-dimensional hypercube, but did not have matching lower bounds.  We
provide lower bounds, showing that their upper bounds are correct to
within constant factors.

\nib
Counting the number of linear extensions of a partially ordered set is
\#P-complete \citep*{brightwell-winkler:extension}, making it
computationally intractible.  (One application of counting linear
extensions is in a data mining application that infers partial orders
\citep*{mannila-meek:partial-order}.)
One can approximately count linear extensions by
randomly generating them, so there have been a number of
articles on generating random linear extensions.  The latest, by
\citet*{bubley-dyer:extension}, proposes a Markov chain in which pairs
of elements in the linear extension are randomly transposed if doing
so respects the partial order.  To make the analysis easy, their Markov
chain selects a random site with probability proportional to parabolic curve,
and then attempts to transpose the elements at this random site.
We show here that the uniform distribution on sites works as well,
obtaining a constant factor that is only marginally worse.

We give further background on these various Markov chains in later 
sections where we analyze them.  In \sref{prelim} we provide basic
definitions, such as what it means formally for the state of a Markov
chain to be close to random.  We study in \sref{path} a Markov
chain for generating random lattice paths, since it is simple, yet
illustrates the key ideas that we will use to analyze the various
other Markov chains.  We use Fourier analysis to define
on the state space of the Markov chain a function $\Phi$ that has a certain
contraction property.  With $S$ denoting the current state of the
Markov chain, and the random variable $S'$ denoting the next state of
the chain, we have $E[\Phi(S')|S]=(1-\gap)\Phi(S)$.  We derive both
upper bounds and lower bounds using this contraction property.  After
\sref{path} the remaining sections may be read in any order.  We see
in \sref{cards} how to apply the results about the path Markov chain
to the chain for shuffling by random adjacent transpositions.  In
\sref{lozenge} we generalize the upper bound for the path Markov chain
to upper bound the mixing time of the lozenge-tiling Markov chain
introduced by \lrs.  In \sref{k&k} we modify \bd's path-coupling
analysis of the \kdk Markov chain to obtain the $O(n^3\log n)$ mixing
time bound.  When using a local randomizing operation to update a high
dimensional configuration, one typically either updates a random
co\"ordinate each step, or else updates the co\"ordinates in a
systematic order.  In \sref{sweep} we compare these two methods for
the chains that we are studying; our analysis indicates that the second method
is better.  We take a second look at the lattice path and
permutation Markov chains in \sref{path-perm}, and refine our previous
arguments to obtain tighter constants.  We consider exclusion and
exchange processes in \sref{exclusion}, where among other things we
resolve the aforementioned questions of \dsc.  Many of the mixing time
upper and lower bounds we give differ by small constant factors.  We
give in \sref{heuristic} heuristic arguments and present experimental
evidence for determining the correct constant factors in the mixing
times.  We summarize in \tbref{times} many of these mixing time bounds
and their (conjectural) correct values.  \sref{heuristic} also contains
several open problems for further research.  We make some concluding
remarks in \sref{conclude}.

\newpage

\newcommand{\logn}{\log n}%{(\log n + c)}
\newcommand{\logl}{\log \ell}%{(\log \ell + c)}
\newcommand{\withinconstants}{\begin{tabular}{c}sometimes tight to\\[-5pt] within constants\end{tabular}}
\begin{table}[ht]
\label{tbl:times}
\begin{tabular}{|r|c||c|c|}
\hline
\multicolumn{4}{|c|}{State space, Markov chain}\\
\hline
  parameter &  \begin{tabular}{@{}c@{}}(conjectural) \\correct answer\end{tabular} & \begin{tabular}{@{}c@{}}rigorous\\lower bound\end{tabular} & \begin{tabular}{@{}c@{}}rigorous\\upper bound\end{tabular} \\
\hline
\hline
\multicolumn{4}{|c|}{Paths in $n/2 \times n/2$ box, random adjacent transpositions} \\
\hline
&&&\\[-11pt]
  coupling threshold &  $2/\pi^2 n^3 \logn$ & $2/\pi^2 n^3 \logn$ & $2/\pi^2 n^3 \logn$ \\
  variation threshold&  $1/\pi^2 n^3 \logn$ & $1/\pi^2 n^3 \logn$ & $2/\pi^2 n^3 \logn$ \\
  separation threshold&  $2/\pi^2 n^3 \logn$ & $1/\pi^2 n^3 \logn$ & $4/\pi^2 n^3 \logn$ \\
  spectral gap       &  $\displaystyle\frac{1-\cos(\pi/n)}{n-1}$
                     &  $\displaystyle\frac{1-\cos(\pi/n)}{n-1}$
                     &  $\displaystyle\frac{1-\cos(\pi/n)}{n-1}$ \\[8pt]
\hline
\hline
\multicolumn{4}{|c|}{Permutations --- $S_n$, random adjacent transpositions} \\
\hline
&&&\\[-11pt]
  coupling threshold &  $4/\pi^2 n^3 \logn$  & $2/\pi^2 n^3 \logn$ & $4/\pi^2 n^3 \logn$ \\
  variation threshold&  $1/\pi^2 n^3 \logn$  & $1/\pi^2 n^3 \logn$ & $2/\pi^2 n^3 \logn$ \\
  separation threshold&  $2/\pi^2 n^3 \logn$  & $1/\pi^2 n^3 \logn$ & $4/\pi^2 n^3 \logn$ \\
  spectral gap       &  $\displaystyle\frac{1-\cos(\pi/n)}{n-1}$
                     &  $\displaystyle\frac{1-\cos(\pi/n)}{n-1}$
                     &  $\displaystyle\frac{1-\cos(\pi/n)}{n-1}$ \\[8pt]
\hline
\hline
\multicolumn{4}{|c|}{Lozenge tilings of order $\ell$ hexagon, Luby-Randall-Sinclair chain} \\
\hline
&&&\\[-11pt]
  coupling threshold & ??  & $8/\pi^2 \ell^4 \logl$ & $48/\pi^2 \ell^4 \logl$ \\
  variation threshold&  $16/\pi^2 \ell^4 \logl$  & $8/\pi^2 \ell^4 \logl$ & $48/\pi^2 \ell^4 \logl$ \\
  separation threshold&  $32/\pi^2 \ell^4 \logl$  & $8/\pi^2 \ell^4 \logl$ & $96/\pi^2 \ell^4 \logl$ \\
  spectral gap    &  $\displaystyle\frac{1-\cos(\pi/(2\ell))}{\ell(2\ell-1)}$
                  &  $\displaystyle\frac{1-\cos(\pi/(2\ell))}{\ell(2\ell-1)}$
                  &  $\displaystyle\frac{1-\cos(\pi/(2\ell))}{\ell(2\ell-1)}$\\[8pt]
\ifodd 1
\hline
\hline
\multicolumn{4}{|c|}{Linear extensions of partially ordered set, Karzanov-Khachiyan chain} \\
\hline
&&&\\[-11pt]
  variation threshold& depends on poset & \withinconstants
%\begin{tabular}{c}$\exists$ poset which needs\\ $1/\pi^2 n^3 \logn$\end{tabular}
 & $4/\pi^2 n^3 \logn$ \\
  spectral gap       &  depends on poset
                     &  $\displaystyle\frac{1-\cos(\pi/n)}{n-1}$
                     &  sometimes tight \\[8pt]
\hline
\hline
\multicolumn{4}{|c|}{Lozenge tilings of region of $n$ triangles and width $w$, Luby-Randall-Sinclair chain} \\
\hline
&&&\\[-11pt]
  variation threshold& depends on region & \withinconstants & $3/\pi^2 w^2 n \logn$ \\
  spectral gap    &  depends on region 
                  &  $\displaystyle\frac{1-\cos(\pi/w)}{n}$
                  &  \withinconstants \\[8pt]
\fi
\hline
\end{tabular}
\caption{Summary of mixing time bounds for several classes of Markov chains considered here.  The variation and separation thresholds are defined in \sref{prelim}.  The bounds for lattice paths and permutations are proved in \sref{path} and \sref{path-perm}, the bounds for lozenge tilings are proved in \sref{lozenge}, and the bounds for the \kdk chain are proved in \sref{k&k}.  The coupling times are for the natural monotone grand couplings described in \sref{path}, \sref{cards}, and \sref{lozenge}.  The conjectural correct answers are derived in \sref{heuristic}.  The spectral gap for permutations was previously known (Diaconis, unpublished).}
\end{table}

\vspace*{-24pt}
\section{Preliminaries}\label{sec:prelim}

Here we review some basic definitions and properties
pertaining to mixing times and couplings.  For a more complete
introduction to these ideas see
\cite{aldous:group-walks} or \cite{aldous-fill:book}.

When a Markov chain is started in a state $x$ and run for $t$
steps, we denote the distribution of the state at time $t$ by $P_x^t$,
where $P$ is the state transition matrix of the Markov chain.  If the
Markov chain is connected and aperiodic, then as $t\rightarrow\infty$
the distribution $P_x^t$ converges to a unique stationary distribution
which is often denoted by $\pi$.  Since we will use $\pi$ to denote
the ratio of the circumference of a circle to its diameter, here we
use $\mu$ to denote the stationary distribution.  In all our examples
the stationary distribution is the uniform distribution, which we
denote with $U$.

To measure the distance between the distributions $P_x^t$ and $\mu$
one usually uses the total variation distance.  With $\mathcal{X}$
denoting the state space, the variation distance is defined by
$$ \| P_x^t - \mu\|_{\TV} = \max_{A\subset\mathcal{X}} |P_x^t(A)-\mu(A)| = \frac12\sum_{y\in\mathcal{X}} |P_x^t(y)-\mu(y)| = \frac12 \| P_x^t - \mu\|_1\ ,$$
The variation distance when the chain is started from the worst start state is denoted by $$d(t) = \max_x \|P_x^t - \mu\|_{\TV}\ ,$$
though it is often more convenient to work with
$$\bar{d}(t) = \max_{x,y} \|P_x^t - P_y^t\|_{\TV}\ $$
since $\bar{d}(t)$ is submultiplicative whereas $d(t)$ is not.  It is easy to see that $d(t)\leq \bar{d}(t)\leq 2 d(t)$.

The (variation) mixing time is the time it takes for $\bar{d}(t)$ to
``become small'' (say less than $1/e$).  It is a surprising fact that
for many classes Markov chains there is a threshold time $T$ such that
$\bar{d}((1-\varepsilon)T)>1-\varepsilon$ but
$\bar{d}((1+\varepsilon)T)<\varepsilon$, where $\varepsilon$ tends to
$0$ as the ``size'' of the Markov chain gets large; see \cite{diaconis:cutoff}
for a survey of this ``cutoff phenomenon.''

Most of our mixing time upper bounds are derived via coupling arguments.
In a (pairwise) coupling
there are two copies $X_t$ and $Y_t$ of the Markov chain that are run
in tandem.  The $X_t$'s by themselves follow the transition rule of
the Markov chain, as do the $Y_t$'s, but the joint distribution of
$(X_{t+1},Y_{t+1})$ given $(X_t,Y_t)$ is often contrived to make the
two copies of the Markov chain quickly coalesce (become equal).  
It is a standard fact that
$$\bar{d}(t) \leq \max_{x,y} \Pr[X_t\neq Y_t|X_0=x,Y_0=y]\ ,$$
so that a coupling which coalesces quickly can give us a good upper
bound on the mixing time.  Most of the variation threshold upper
bounds in \tbref{times} follow from this relation and a corresponding
coupling time bound, and lower bounds on the coupling time likewise
follow from corresponding lower bounds on the variation threshold
time.  The remaining variation threshold upper bounds and coupling
time lower bounds in \tbref{times} that do not follow from this
relation are derived in
\sref{path-perm}.

Many pairwise couplings can be extended to ``grand couplings,'' where
at each time step there is a random function $F_t$ defined on the
whole state space of the Markov chain, and $X_{t+1}=F_t(X_t)$ and
$Y_{t+1}=F_t(Y_t)$ for any $X_t$ and $Y_t$.
For example, if the state space is the set of permutations on $n$ cards,
then the update rule ``pick a random adjacent pair of cards,
and flip a coin to decide whether to place them in
increasing order or decreasing order'' defines a grand coupling;
the choice of the adjacent pair and the value of the coin flip define
the random function on permutations.  In \sref{k&k} and
\sref{path-perm} we will also consider pairwise couplings that do not
extend to grand couplings.

All of the grand couplings considered in this article are monotone,
which is to say that there is a partial order $\preceq$ such that if
$x\preceq y$ then also $F_t(x)\preceq F_t(y)$.  All of the partial
orders considered here have a maximal element, denoted $\1$, and a
minimal element, denoted $\0$, i.e.\ so that $\0\preceq x \preceq \1$
for each $x$ in the state space.  Monotone grand couplings are
particularly convenient for algorithms (see e.g.\
\cite{propp-wilson:exact-sampling} or \cite{fill:interruptible}).

In \sref{heuristic} we consider not only the variation distance,
but also the separation distance, which is defined by $$ s(t) =
\max_{x,y} \frac{\mu(y)-P_x^t(y)}{\mu(y)}\ .$$ The function $s(t)$ is also
submultiplicative, and also often exhibits a sharp threshold.  In
general $d(t)\leq s(t)$, and for reversible Markov chains $s(2t)\leq
2\bar{d}(t)-\bar{d}(t)^2$ (see \cite[Chapt.~4,
Lemma~7]{aldous-fill:book}).  The rigorous bounds in \tbref{times}
pertaining to separation distance follow from these relations and the
corresponding bounds for the variation distance.

\section{Lattice path Markov chain}\label{sec:path}

\begin{figwindow}[0,r,%
{\psfig{figure=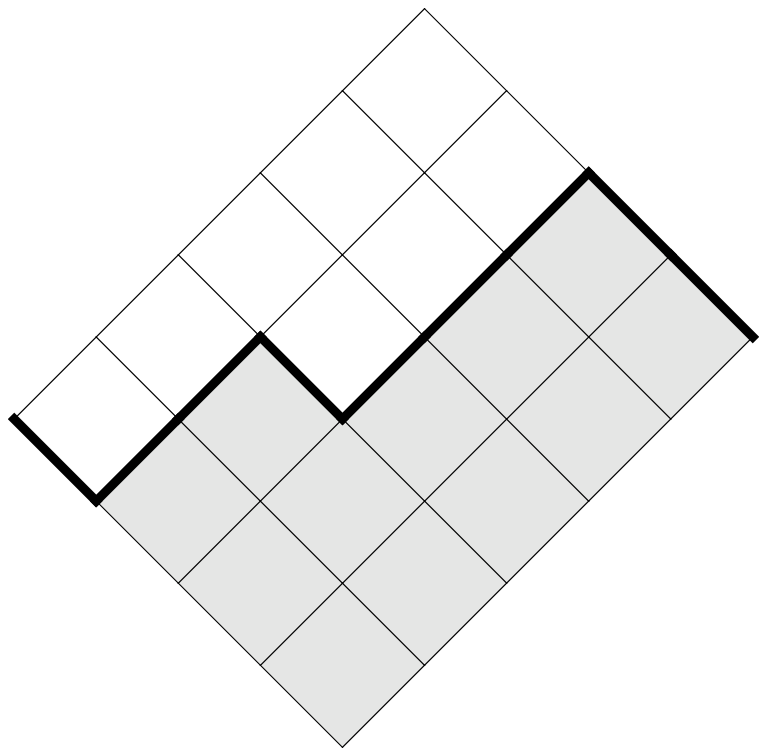}},%
{A lattice path through the $4\times 5$ rectangle, with the $1\times 1$ boxes underneath it shaded.  Here $a=4$, $b=5$, $n=9$, and the lattice path's
encoding is $011011100$.\label{fig:path}}]
\noindent
Consider an $a\times b$ rectangle composed of $1\times 1$ boxes
rotated $45^\circ$, so that the sides of length $a$ are oriented
northwest/southeast.  A lattice path (see \fref{path}) is a traversal
from the left-most corner to the rightmost corner of the rectangle,
traveling along the borders of the $1\times 1$ boxes,
so that each move
is either up-and-right or down-and- right.  Such lattice paths can be
encoded as strings of length $a+b$ consisting of $a$ $0$'s (down moves)
and $b$ $1$'s (up moves).  There are $n!/(a!b!)$ such lat- tice paths,
where for convenience we let $n=a+b$.  These
lattice paths correspond to sets of
$1\times 1$ boxes which are ``stable under gravity,''
i.e.\ no box lies above an empty cell.

\mbox{\hspace{\parindent}}
Consider the following Markov
chain for ran-
domly generating a lattice path between the oppo-
site
corners of the $a\times b$ rectangle.  Given a path, the Markov chain
randomly picks one of the $(n-1)$ internal columns (we assume $n\geq
2$), and then ran-
domly decides whether to try pushing the path up at
that point, or to try pushing it down.  If pushing the path up (or
down) would result in an invalid path, the Markov chain simply sits
idle during that step.  Equivalently, in the binary string representation,
the Markov chain picks a random
adjacent pair of letters in the string, and then randomly decides to
either sort them or reverse-sort them.  Understanding this Markov
chain will be instrumental to understanding the Markov chains for
random lozenge tilings, for card shuffling by random adjacent
transpositions, and for other types of card shuffling.
%  The
%mixing time bound for this example given in
%\citep*{luby-randall-sinclair:markov-lattice} (using
%\lref{diffuse}) is $O(\ell N^2)$, where $N=\ell^2$ is the area, or
%$O(\ell^5)$.  The bound that we get here is $O(n^3 \log \min(a,b))$, which is
%within a constant factor of optimal when $a$ and $b$ are polynomially related.
\end{figwindow}

\subsection{Contraction property}

We will analyze this Markov chain by measuring 1) the ``displacement''
of a given single path on the $a\times b$ square from ``equilibrium'',
and 2) the ``gap'' between two such paths when one of the paths is
entirely above the other.  We lower bound the mixing time by computing
the displacement of a path when not enough steps have been taken, and
showing that it is typically different from the displacement of a
random path.  We upper bound the mixing time by showing that after
enough steps, starting from the top and bottom paths, the expected gap
is so small that they have almost surely coalesced to the same path.

It will be useful later to use horizontal co\"ordinates that range
from $-n/2$ to $n/2$.
Let $h(x)$ denote the height
of the path at position $x$ relative to the line connecting the
opposite corners of the box.  That is, $h(x)$ is the number of up
moves to the left of position $x$ minus the expected number of such up
moves.  Thus $h(-n/2)=h(n/2)=0$, and $h(x) = h(x-1) + a/n$ if there
was an up move between $x-1$ and $x$, or else
$h(x)=h(x-1) - b/n$ if there was a down move.
For the example in
\fref{path}, the heights change by $+4/9=a/(a+b)$ for up moves and $-5/9=-b/(a+b)$ for down moves, and are as follows:\\
\centerline{\begin{tabular}{|c|cccccccccc|}\hline
$x$ &
-9/2&
-7/2&
-5/2&
-3/2&
-1/2&
1/2 &
3/2 &
5/2 &
7/2 &
9/2
\\
\hline
$h(x)$ &
0&
-5/9&
-1/9&
3/9&
-2/9&
2/9&
6/9&
10/9&
5/9&
0\\
\hline \end{tabular}.}

\vspace*{8pt}
The displacement function $\Phi$ of $h$ that we
will find useful is
\begin{equation}\label{eqn:disp}
\Phi(h) = \sum_{x=-n/2}^{n/2} h(x) \cos\frac{\beta x}{n}\ ,
\end{equation}
where $0\leq\beta\leq \pi$.
This function weighs deviations from expectation more heavily
near the middle of the path than near its endpoints.
  Given two lattice
paths with height functions \hl and \hu, where $\hl(x)\leq \hu(x)$ for all
$x$, define the gap function to be $\hu-\hl$, and the gap to be
\begin{equation*}
\Phi(\hl,\hu) = \Phi(\hu-\hl) = \Phi(\hu)-\Phi(\hl)\ .
\end{equation*}
Note that since $0\leq\beta\leq\pi$, the gap is strictly positive when
the paths \hl and \hu differ, and is $0$ otherwise.  After the Markov
chain has equilibrated, so that each path is equally likely,
$E[h(x)]=0$, so the expected displacement is $E[\Phi(h)]=0$.

\begin{lemma}
\label{lem:path-df}
Let the displacement function $\Phi$ be defined by \eref{disp}.  Suppose
$h$ is a height function (so $\Phi(h)$ is the displacement) and
$\beta=\pi$, or $h=\hu-\hl$ is a gap function (so $\Phi(h)$ is the gap) and $0\leq\beta\leq\pi$.  Let $h'$ be the height or gap function after one step of the Markov chain.  Then
$$E[\Phi(h')-\Phi(h)|h] \leq \frac{-1+\cos(\beta/ n)}{n-1} \Phi(h)\ ,$$
with equality when $\beta=\pi$.  The coefficient on the right-hand
     side is bounded by
$$-\frac{\beta^2}{2 n^2(n-1)} \leq  \frac{-1+\cos(\beta/ n)}{n-1} \leq
  -\frac{\beta^2}{2 n^3}\ .$$
\end{lemma}

\begin{proof}
Suppose we pick a site $x$, flip a coin, and adjust the height
accordingly.  Then the expected value of the new height at $x$ is just
$[h(x+1)+h(x-1)]/2$.  Assume that we pick each site (other than $-n/2$
and $n/2$) with probability $1/p$, where $p=n-1$ is the number of positions that can be picked.  Then with primes denoting the
updated variables, $$E[h'(x)|h] = \frac{p-1}{p} h(x) +
\frac{1}{p}\frac{h(x+1)+h(x-1)}{2}$$ when $-n/2 < x < n/2$, so that
\begin{align}
    E[\Phi(h')|h] &= \sum_{x=-n/2+1}^{n/2-1} E[h'(x)|h] \cos\frac{\beta x }{n}
    \label{eqn:phip}\\
    E[\Phi(h')|h] &= \frac{p-1 }{p}\Phi(h) + \frac{1/2}{p}
    \sum_{x=-n/2+1}^{n/2-1} [h(x+1)+h(x-1)] \cos\frac{\beta x }{n} \label{eqn:exp}\\
E[\Phi(h')-\Phi(h)|h]&= \frac{-1 }{p}\Phi(h) + \frac{1/2}{p}
    \sum_{\substack{-n/2+1 \leq x \leq n/2-1\\
                   -n/2 \leq y \leq n/2\\
                   |x-y|=1
         }}             h(y) \cos\frac{\beta x }{n} \notag \\
         &= \frac{-1 }{p}\Phi(h) + \frac{1/2}{p}
    \sum_{\substack{-n/2+1 \leq x \leq n/2-1\\
                   -n/2+1 \leq y \leq n/2-1\\
                   |x-y|=1
         }}             h(y) \cos\frac{\beta x }{n} \label{eqn:h0}\\
         &\leq \frac{-1 }{p}\Phi(h) + \frac{1/2}{p}
    \sum_{\substack{-n/2   \leq x \leq n/2\\
                   -n/2+1 \leq y \leq n/2-1\\
                   |x-y|=1
         }}             h(y) \cos\frac{\beta x }{n} \label{eqn:cos}\\
            &= \frac{-1 }{p}\Phi(h) + \frac{1/2}{p}
    \sum_{y=-n/2+1}^{n/2-1} h(y) \left[\cos\frac{\beta(y-1)}{n} + \cos\frac{\beta(y+1)}{n}\right] \notag\\
\intertext{}
            &= \frac{-1 }{p}\Phi(h) + \frac{1/2}{p}
    \sum_{y=-n/2+1}^{n/2-1} h(y) 2\cos\frac{\beta y}{n} \cos\frac{\beta}{n} \label{eqn:trig}\\
            &= \frac{-1 + \cos(\beta/n) }{p} \Phi(h) = \frac{-1 + \cos(\beta/n) }{n-1} \Phi(h) \notag\ ,
\end{align}
where we have used $E[h'(\pm n/2)]=0$ in Equation~\eqref{eqn:phip},
$h(\pm n/2)=0$ in Equation~\eqref{eqn:h0}, and the trigonometric
identity $\cos(\theta-\phi) + \cos(\theta+\phi) = 2\cos(\theta)\cos(\phi)$ in
Equation~\eqref{eqn:trig}.  Inequality~\eqref{eqn:cos} is justified if
$\beta\leq\pi$, since by assumption $h=\hu-\hl\geq 0$ for each
$x$, and $\cos(\beta/2)\geq 0$; when $\beta=\pi$ it becomes an equality.

To upper bound the right-hand side we use
$\cos(x) \leq 1-x^2/2+x^4/24 = 1+(x^2/2)(-1+x^2/12)$.
\begin{align*}
\frac{-1 + \cos(\beta/n) }{n-1}
&\leq \frac{1}{n-1} \frac{(\beta/ n)^2}{2} \left(-1+\frac{(\beta/n)^2}{12} \right) \\
&\leq -\frac{\beta^2}{2 n^3}\ ,
\end{align*}
where to get the last line we used
$(1-c/n^2)/(n-1) \geq 1/n$ when $n> 1$ and $-c/n\geq -1$.
Here $c\leq \pi^2/12$, so these conditions are satisfied whenever $n> 1$.
The lower bound is somewhat easier: use the bound
$\cos(x)\geq 1-x^2/2$.
\end{proof}

\subsection{Upper bound}

\begin{theorem} \label{thm:path-upper}
When $n$ is large, after $$\frac{2+o(1)}{\pi^2} n^3
\log\frac{ab}{\varepsilon}$$ steps the variation distance
from stationarity is $\varepsilon$, and the probability that the two extremal
paths have coalesced is $1-\varepsilon$.  (The $o(1)$ term is a
function of $n$ alone.)
\end{theorem}

\citet*{felsner-wernisch:linear-extension} reference an early
version of this paper, since they used the upper bound in this theorem
to bound the rate of convergence of the \kdk Markov chain for
generating random linear extensions of a certain class of partially
ordered sets.  The \kdk Markov chain is discussed further in \sref{k&k}.

\begin{proof}[Proof of \tref{path-upper}]
To obtain the upper bound, we consider a pair of coupled paths $\hu_t$
and $\hl_t$ such that $\hu_0$ is the topmost path, $\hl_0$ is the
bottommost path.  The sequences $\hu_t$ and $\hl_t$ are generated by
the Markov chain using ``the same random moves'', so that $\hu_{t+1}$
and $\hl_{t+1}$ are obtained from $\hu_t$ and $\hl_t$ respectively by
sorting or unsorting (same random decision made in both cases) at the
same random location $x$.  One can check by induction that $\hu_t\geq\hl_t$.
Let $\Phi_t = \Phi(\hu_t-\hl_t)$; $\Phi_t=0$ if and only if $\hu_t=\hl_t$.

We compute $E[\Phi_t]$; when it is small compared to the minimum
possible positive value of $\Phi_t$, it will follow that with high
probability $\Phi_t=0$.  By choosing $\beta$ to be slightly smaller
than $\pi$, we make this minimum positive value not too small, and
thereby get a somewhat improved upper bound.

From \lref{path-df} and induction, we get
\begin{equation*}
E[\ft] \leq \Phi_0 \left[1 - \frac{1-\cos(\beta/ n) }{n-1}\right]^t
          \leq \Phi_0 \exp \left[- \frac{\beta^2 t }{2 n^3} \right]\ .
\end{equation*}
But $E[\ft] \geq \Pr[\ft>0] \fm$.
Thus after $t\geq (2/\beta^2) n^3\log (\Phi_0/(\fm \varepsilon))$ steps
$\Pr[\ft>0] \leq \varepsilon$.
We have $\Phi_0 \leq a b$, and $\fm=\cos(\beta (n/2-1)/n) >
\cos(\beta/2) \approx (\pi-\beta)/2$.  The optimal choice of $\beta$
is $\pi - \Theta(1/\log n)$, but all that matters is that
$\pi-\beta\rightarrow 0$ as $n\rightarrow\infty$ while $\log
(1/(\pi-\beta)) \ll \log (ab)$.  Substituting, we find that $t =
(2/\pi^2+\Theta(\log\log n/\log n)) n^3\log (ab/\varepsilon)$ steps are
enough to ensure that the probability of coalescence is at least
$1-\varepsilon$.
\end{proof}

\remark
We will show in \sref{path-perm} that when $a=b=n/2$
the coupling time is actually
$(2/\pi^2) n^3\log n$.

\remark
Proving mixing time upper bounds via a contraction property in the ``distance'' between configurations is a fairly standard technique.  Traditionally the distance has been measured in terms of Hamming distance or other integer-valued distance, which for our applications does not yield the requisite contraction property.

We get the spectral gap entries in \tbref{times} using similar reasoning
(see also, e.g., \cite{chen:coupling-eigenvalue}):
\begin{proposition}
  If a function $\Phi$ is strictly monotone increasing in the partial
  order of a reversible monotone Markov chain with top and
  bottom state, and whenever $X\preceq Y$ we have
  $E[\Phi(Y')-\Phi(X')|X,Y]\leq(1-\gamma)(\Phi(Y)-\Phi(X))$, then the
  spectral gap must be at least $\gamma$.
\end{proposition}
\begin{proof}
  Perturb the stationary distribution by an eigenvector associated
  with the second largest eigenvalue $\lambda$, and run the Markov
  chain starting from this distribution.  After $t$ steps the
  variation distance from stationarity is $A \lambda^t$.  If we run
  the Markov chain starting from the top and bottom states, after
  $t$ steps the states are different with probability at most
  $(\Phi(\1)-\Phi(\0))(1-\gamma)^t/\min_{X\prec Y}(\Phi(Y)-\Phi(X))$.
  The coupling time bound on the variation distance gives
  $\lambda\leq 1-\gamma$.
\end{proof}

\subsection{Lower bound}

We obtain a lower bound on the mixing time when the rectangle is not
too narrow:

\begin{theorem} \label{thm:path-lower}
If $\min(a,b)\gg 1$, then after 
$$\frac{1-o(1)}{\pi^2} n^3\log \min(a,b)$$
steps the variation distance from stationarity is $1-o(1)$.
\end{theorem}

We use \lref{path-df} with $\beta=\pi$ so that we get an exact
expression for $E[\Phi_t]$.  Then we bound $\Var[\Phi_t]$ to show that
the distribution of $\Phi_t$ is sharply concentrated about its
expected value.  When $\Phi_t$ and $\Phi_\infty$ are sharply
concentrated about values that are far enough apart, the chain is far
from equilibrium.
(This second-moment approach was also used by
\citet*{diaconis-shahshahani:cube-convergence} to lower
bound the mixing time of random walk on $\Z_2^d$, and by
\citet*{lee-yau:log-sobolev} to lower bound the mixing time of the exclusion
process on a circle.)  The following technical lemma formalizes this argument,
and is used to derive the mixing time lower bounds we give in this article,
except the bound proved in \sref{dsc}, where we need a generalization.
(The coupling time lower bound proved in \sref{path-perm-lower} uses a
different approach altogether.)

\begin{lemma}\label{lem:anticonverge}
If a function $\Phi$ on the state space of a Markov chain satisfies $E[\Phi(X_{t+1})|X_t] = (1-\gap) \Phi(X_t)$, and $E[(\df)^2|X_t]\leq R$ where $\df=\Phi(X_{t+1})-\Phi(X_t)$, then when the number of Markov chain steps
$t$ is bounded by
$$ t \leq \frac{\log \fx + \frac{1}{2} \log \frac{\gap \varepsilon }{4 R}
  }{- \log(1-\gap)}\ ,$$
and $0 < \gap \leq 2-\sqrt{2} \doteq 0.58$ (or else $0<\gap\leq 1$
  and $t$ is odd),
then the variation distance from stationarity is at least $1-\varepsilon$.
\end{lemma}

Before proving \lref{anticonverge}, we show how to use it to prove \tref{path-lower}.

\begin{proof}[Proof of \tref{path-lower}]
By \lref{path-df}, our function $\Phi$ satisfies the contraction
property required by \lref{anticonverge} when $\beta=\pi$, and $$
\gap = \frac{1-\cos(\pi/ n) }{n-1} \mathrel{\overset{\geq}{\underset{\sim}{\text{and}}}}
\frac{\pi^2}{2 n^3}\ .$$ The constraint on $\gap$ is satisfied when
$n\geq 3$.

To get a bound $R$, observe that any path $h$ can have at most
$2\min(a,b)$ local extrema, so $\Pr[\df\neq 0] \leq \min(a,b)/(n-1)$.
But $|\df|\leq 1$, so $\max_h E[(\df(h))^2] \leq \min(a,b)/(n-1) \equiv R$.

The maximal path maximizes $\Phi$, giving $\Phi_0 = \Theta(a b)$.
Substituting into \lref{anticonverge} and simplifying, $\gap/R \sim
1/(n^2\min(a,b))$, so the numerator becomes $\log(a b) - \log n -
(1/2)\log \min(a,b)+O(1) = \log\min(a,b)+O(1)$ for bounded values of
$\varepsilon$, giving our lower bound of $(1/\pi^2-o(1))n^3\log
\min(a,b)$.
\end{proof}

\begin{proof}[Proof of \lref{anticonverge}]
Let $\Phi_t = \Phi(X_t)$.
By induction $$E[\Phi_t)|X_0] = \Phi_0 (1-\gap)^t\ .$$
By our assumptions on $\gap$, in equilibrium $E[\Phi]=0$.

With $\df$ denoting $\fp-\ft$, we have
\begin{align*}
\fp^2       &= \ft^2 + 2\ft\df + (\df)^2\\
E[\fp^2|\ft]&= (1-2\gap) \ft^2 + E[(\df)^2|\ft] \leq (1-2\gap)
            \ft^2 + R\ ,
\intertext{and so by induction,}
E[\ft^2]    &\leq \Phi_0^2 (1-2\gap)^t + \frac{R}{2\gap}\ ,
\intertext{then subtracting $E[\ft]^2$,}
\Var[\ft]   &\leq \Phi_0^2 [(1-2\gap)^t-(1-\gap)^{2t}] 
                  + \frac{R}{2\gap} \\
\Var[\ft]   &\leq \frac{R}{2\gap}
\end{align*}
for each $t$.  To get the last line we used our constraints on
$\gap$ and $t$: $(1-\gap)^2 = 1-2\gap+\gap^2 \geq
1-2\gap$, so when $t$ is odd, $(1-\gap)^{2t} \geq
(1-2\gap)^t$.  When $t$ is even, we need $(1-\gap)^2 \geq
2\gap-1$ as well, which is satisfied when $\gap \leq 2-\sqrt{2}$
or $\gap \geq 2+\sqrt{2}$.

From Chebychev's inequality, $$\Pr\left[|\ft-E[\ft]| \geq
  \sqrt{R/(2\gap\varepsilon)}\right] \leq \varepsilon\ .$$
As $E[\Phi_\infty]=0$, if $E[\ft] \geq \sqrt{4 R/(\gap\varepsilon)}$,
then the probability that $\ft$ deviates below
$\sqrt{R/(\gap\varepsilon)}$ is at most $\varepsilon/2$, and the
probability that $\Phi$ in stationarity deviates above this threshold
is at most $\varepsilon/2$, so the variation distance between the
distribution at time $t$ and stationarity must be at least
$1-\varepsilon$.
If we take the initial state to be the one maximizing $\Phi_0$,
then $$E[\ft] = \fx (1-\gap)^t \geq \sqrt{4 R/(\gap\varepsilon)}$$
when
\begin{align*}
 t &\leq \frac{\log \left[\fx\div\sqrt{\frac{4 R}{\gap\varepsilon}}\right]}{-\log(1-\gap)}\ .\qedhere
\end{align*}
\end{proof}

\subsection{Intuition}
Since the expected value of the new height function is a certain local
average of the current height function, the evolution of the height
function $h$ (or rather its expected value) proceeds approximately
according to the rule
\begin{equation}\label{eq:heat}
\frac{\partial h }{\partial t} = \text{const.}
\times \frac{\partial^2 h }{\partial x^2}\ .
\end{equation}
Since the equation is linear in $h$, it is natural to consider its
eigenfunctions, which are just the sinusoidal functions that are zero
at the boundaries.  We can decompose any given height function into a
linear combination of these sinusoidal components, and consider the
evolution of each component independently.  The displacement function
$\Phi$ (when $\beta=\pi$) is just the coefficient of the principal
mode (the sinusoidal function with longest period) when we decompose
the height function in this way.  The coefficient of the principal
mode decays the most slowly, making it the most useful for purposes of
establishing a lower bound.  For purposes of establishing an upper
bound, the difference between two height functions also obeys
Equation~\eqref{eq:heat}.
We used the fact that when two paths, one above the other,
have the same displacement, then they are the same path.  The
coefficient of the principal mode is the only one for which we can
guarantee this property, since the other eigenfunctions take on both
positive and negative values in the interior.  Thus the principal mode
essentially controls the rate of convergence, making its coefficient a
natural displacement function.

\section{Card shuffling by adjacent transpositions}\label{sec:cards}

\begin{figwindow}[0,r,%
{\psfig{figure=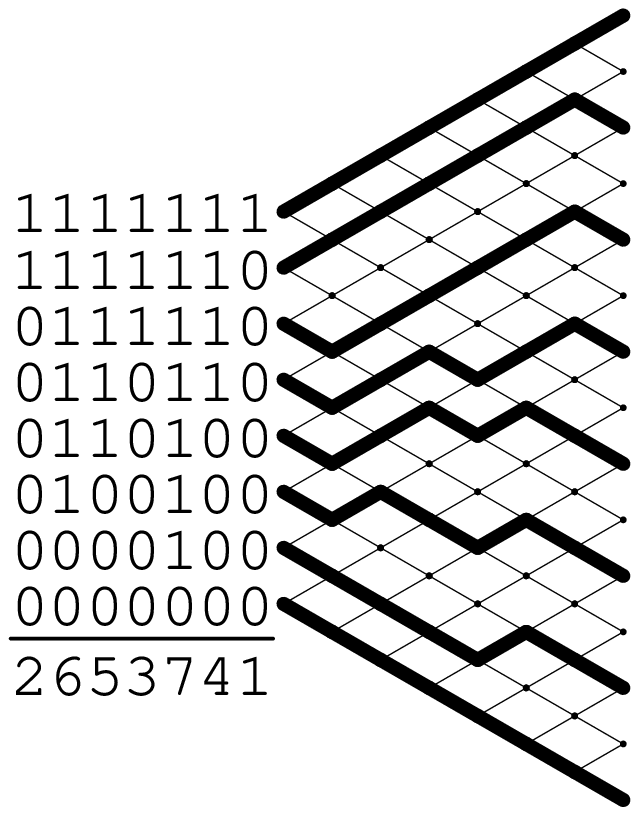}},%
{The permutation $2653741$ and its $8$ threshold functions shown as lattice paths.\label{fig:permutation}}]
\noindent
Next, we analyze the card shuffling Markov chain that trans-
poses random
adjacent pairs of cards.  (This Markov chain is a special case of 
the move-ahead-one update rule that has been studied in self-organizing
linear search; \citet{hester-hirschberg:self-organizing} give a survey.)
We will consider this Markov chain to be
implemented according to the rule: pick a random adjacent pair
$(i,i+1)$ of cards, flip a coin $c$, and then sort the items in that
adjacent pair if heads, other-
wise reverse-sort them.  The same random
update defined by $i$ and $c$ may be applied to more than one
permutation to obtain coupled Markov chains.  A permutation on the
num-
bers $1,\ldots,n$ has associated with it $n+1$ threshold functions,
where the $i$th threshold function ($0\leq i\leq n$) is a string of
$i$ 1's and $n-i$ 0's, with the 1's at the locations of the $i$
largest numbers of the permutation.  The permutation can be recov-
ered
from these threshold functions simply by adding them up (see
\fref{permutation}).  When a random adjacent pair of numbers within
the permutation are transposed, the effect on any given threshold
function is to transpose the same adjacent pair of 0's and 1's.  The
identity permutation $12\cdots n$ and its reverse $n(n-1)\cdots 1$
give the minimal and maximal paths for any threshold function.  So
when the coupled Markov chains started at these two permutations
coalesce (take on the same value), the grand coupling would take
any starting permutation to this same value.
We can therefore use our analysis of the Markov chain on lattice paths
to analyze the Markov chain on permutations.
\end{figwindow}

\begin{theorem}
\label{thm:cards}
After $(1/\pi^2-o(1)) n^3\log n$ shuffles, the variation distance from
stationarity is $1-o(1)$, and the probability of coalescence is
$o(1)$.  After $(6/\pi^2+o(1)) n^3\log n$ shuffles, the variation
distance from stationarity is $o(1)$, and the probability of
coalescence is $1-o(1)$.
\end{theorem}

We will prove better upper bounds \sref{k&k} and \sref{path-perm-upper};
the point here is to give a quick and easy proof.

\begin{proof}[Proof of \tref{cards}]
The lower bound comes from considering the $\lfloor n/2 \rfloor$th
threshold function.  By Thereom~\ref{thm:path-lower}, after $[1/
\pi^2-o(1)] n^3\log n$ steps the variation distance of just this one
threshold function from stationarity is $1-o(1)$.  The variation
distance from stationarity of the permutation itself is at least as
large.

\hspace{\parindent} The upper bound follows from \tref{path-upper} when we take
$\varepsilon = \delta/n$.  As $ab/(\delta/n) \leq n^3/\delta$, after
$[2/\pi^2 + o(1)] n^3 \log (n^3/\delta)$ steps the probability of any
one given threshold function differing for the upper and lower
permutations is $\leq \delta/n$.  The probability that the upper and
lower permutations differ is at most the expected number of threshold
functions for which they differ, which is at most $\delta$.  Taking
$\delta \ll 1$ but $\log (1/\delta) \ll \log n$, after $[6/\pi^2 + o(1)]
n^3 \log n$ steps coalescence occurs with probability $1-o(1)$.
\end{proof}

\section{Lozenge tilings}\label{sec:lozenge}

\subsection{Background}

Random tilings, or equivalently random perfect matchings, were
originally studied in the physics community as a model for dimer
systems (see e.g.\ \citep{fisher:dimer} and \citep{kasteleyn:phase} and
references contained therein).  Physicists were interested in properties
such as long-range correlations in the dimers.  In the case of
lozenge tilings, the dimers are the lozenges, and the monomers are the
two regular triangles contained in a lozenge.  In recent years
mathematicians have been studying random lozenge tilings (among other
types of random tilings), and have proved a number of difficult
theorems about their asymptotic properties.  For instance, when a very
large hexagonal region is randomly tiled by lozenges, with high
probability the tiling will exhibit a certain circular shape, and the
density of each of the three orientations of lozenges, as a function
of position, is also known \citep*{cohn-larsen-propp:hexagon}.
(See also \citep*{cohn-kenyon-propp:variational}.)
Observations of random lozenge tilings of very large regions played an
important role in the history of these theorems, since they indicated
what sort of results might be true before they were proved, thereby
guiding researchers in their efforts.

Consequently, there have been several articles
\citep*{propp-wilson:exact-sampling}
\citep*{luby-randall-sinclair:markov-lattice}
\citep*{wilson:det}
\citep*{ciucu-propp:shuffling} on techniques to randomly generate
lozenge tilings and other types of tilings.  The first two of these articles
use a Markov chain approach, while the second two use linear algebra.
The first of these articles 
\citep*{propp-wilson:exact-sampling}
introduces monotone-CFTP, which lets one efficiently generate random
structures (e.g.\ lozenge tilings) using special Markov chains,
without requiring any knowledge about the convergence rate of the
Markov chain.  It is the article by
\citet*{luby-randall-sinclair:markov-lattice}
that is most relevant to us here.  In it they introduce novel Markov
chains for generating lozenge tilings (and two other types of
structures).  In this case knowledge of the mixing time of these
Markov chains does not help with the specific task of random generation, as
monotone-CFTP, which determines on its own how long to run a Markov
chain, may be applied to each of these
Markov chains.  But there are still several reasons to determine
the mixing time: (1) in the same way that designers of efficient
algorithms like to prove that the algorithms actually are efficient,
it is desirable to have a proof that the Markov chain is rapidly mixing,
(2) there are physical interpretations of the mixing properties of
dimer systems (see the discussion of \cite{destainville:rhombus} and
references contained therein), and (3)
there has been some speculation \citep{propp:personal} that knowledge
of the convergence properties of these Markov chains may be converted
into knowledge about random tilings of the whole plane (but this
remains to be seen).  For these reasons, \lrs establish
polynomial time bounds on the convergence rates of each of their
Markov chains.  In this section we substantially improve the analysis
of the lozenge tiling Markov chain, and in many cases our bounds
differ by just a constant factor from the true convergence rate.
But as discussed in \sref{other-tower}, \lrs also analyzed other
Markov chains for which it is not clear how to
apply the methods of this section.
Nonetheless, empirical studies suggest that these other Markov
chains converge about as quickly as the lozenge tiling Markov chain.

There is a well-known correspondence between dimers on the hexagonal
lattice, lozenge tilings, and nonintersecting lattice paths of the
type we considered in \sref{path}.  \fref{bijections} illustrates this
correspondence by showing a random perfect matching of a region of the
hexagonal lattice, an equivalent random tiling of a related region by
lozenges, and an equivalent random collection of nonintersecting
lattice paths.  Following \lrs we shall use the lattice path representation of
lozenge tilings.

\begin{figure}[htbp]
\centerline{\psfig{figure=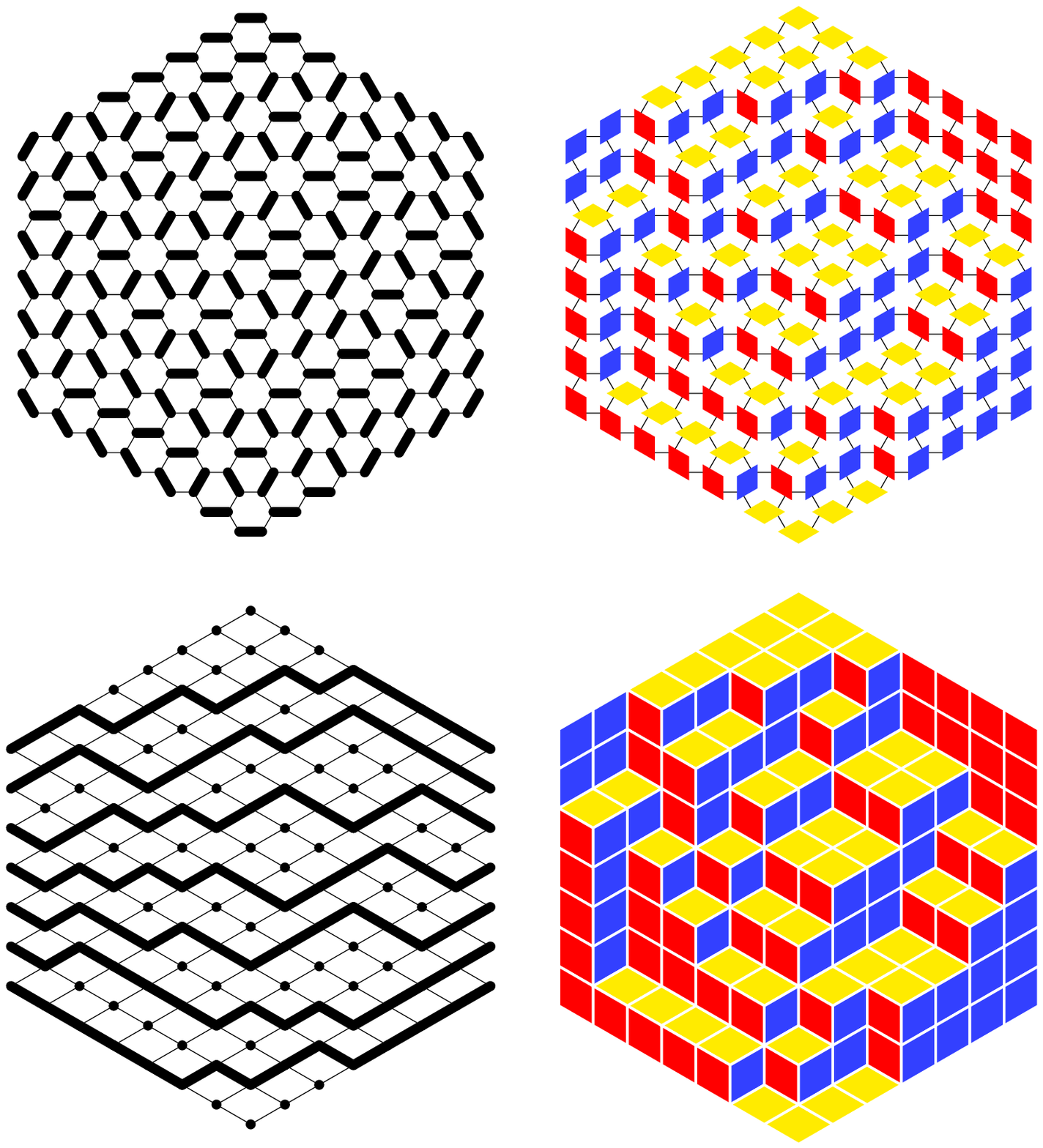}}
\caption{Shown starting in the upper left and proceeding clockwise are
1) a random perfect matching (every vertex paired with exactly one
neighbor) of a portion of the hexagonal lattice,
2) the same perfect matching where edges of the matching shown as
small lozenges, 3) same as 2, but the lozenges are large enough to
touch one another, forming a lozenge tiling of a certain region, and
4) same as 3, but horizontal lozenges are represented as dots, while
the other types are represented as ascending or descending line
segments, which form nonintersecting lattice paths.  These
transformations are bijective, so that any set of nonintersecting
lattice paths corresponds to a lozenge tiling which in turn
corresponds to a perfect matching of the original hexagonal lattice graph.}
\label{fig:bijections}
\end{figure}

\subsection{Displacement function}

In this section we apply the same techniques used to upper bound the
mixing time of the lattice path Markov chain to upper bound the mixing
time of a Markov chain for generating random lozenge tilings.  There
are several Markov chains for generating random lozenge tilings of
regions (each of which possesses the monotonicity property required by
monotone-CFTP) --- the one that we shall analyze was introduced by
\lrs.  Two of these Markov chains use the lattice path representation
of lozenge tilings, and may be viewed as generalizations of the
lattice path Markov chain that we studied already.

Consider the Markov chain that picks a random point on some lattice
path, and then randomly decides whether to try pushing it up or down.
The Markov chain is connected because the top path can be pushed to
its maximum height, then the next highest path, and so on, so that
each configuration can reach a unique maximal configuration.
(Similarly, there is a unique minimal configuration.)  The Markov
chain is aperiodic since pushing the same lattice point twice in the
same direction results in no change the second time.  The Markov
chain is symmetric, so its unique stationary distribution is the
uniform distribution.

\remark
\lrs assumed that the region to be tiled is simply connected, since
otherwise the lattice paths cannot cross the interior holes in the
region, causing the state space to be disconnected.  But if we
restrict the state space to configurations with a specified set of
lattice paths passing under each interior hole, then this restricted
state space is connected by these local moves.  Our mixing time upper
bounds will apply to each such connected component of the state space
if the region to be tiled is not simply connected.

Such is the ``local moves'' Markov chain for lozenge tilings.
Unfortunately, it is difficult to analyze in the same way that we
analyzed the path Markov chain because the paths must remain disjoint.
We would like to define the displacement $\Phi$ of the lozenge tiling
to be the sum of the displacements of each of its paths, but computing
$E[\df]$ is difficult.  We cannot compute the expected new height of a
point on a lattice path simply by looking at the heights of its
neighbors, since another lattice path may or may not be nearby and
block its movement.  \lrs introduced ``nonlocal moves'' to circumvent
this problem and make a Markov chain that is more tractable.

Consider a single lattice path site in isolation.  Its height will
change only if it is a local extremum and gets pushed up (if it's a
minimum) or down (if it's a maximum).  The idea behind the nonlocal
moves is to preserve the expected change in height even if there are
other lattice paths that might block the movement of this local
extremum.  If there are $k$ paths blocking the movement of the local
extremum, then with probability $1/(k+1)$ the corresponding points in
the $k+1$ paths are each moved.  (Naturally if the border of the
region itself prevents the movement of the paths, then the probability
remains zero.)  See \fref{tower}.  This modified Markov chain is
\begin{figure}[htbp]
\centerline{\psfig{figure=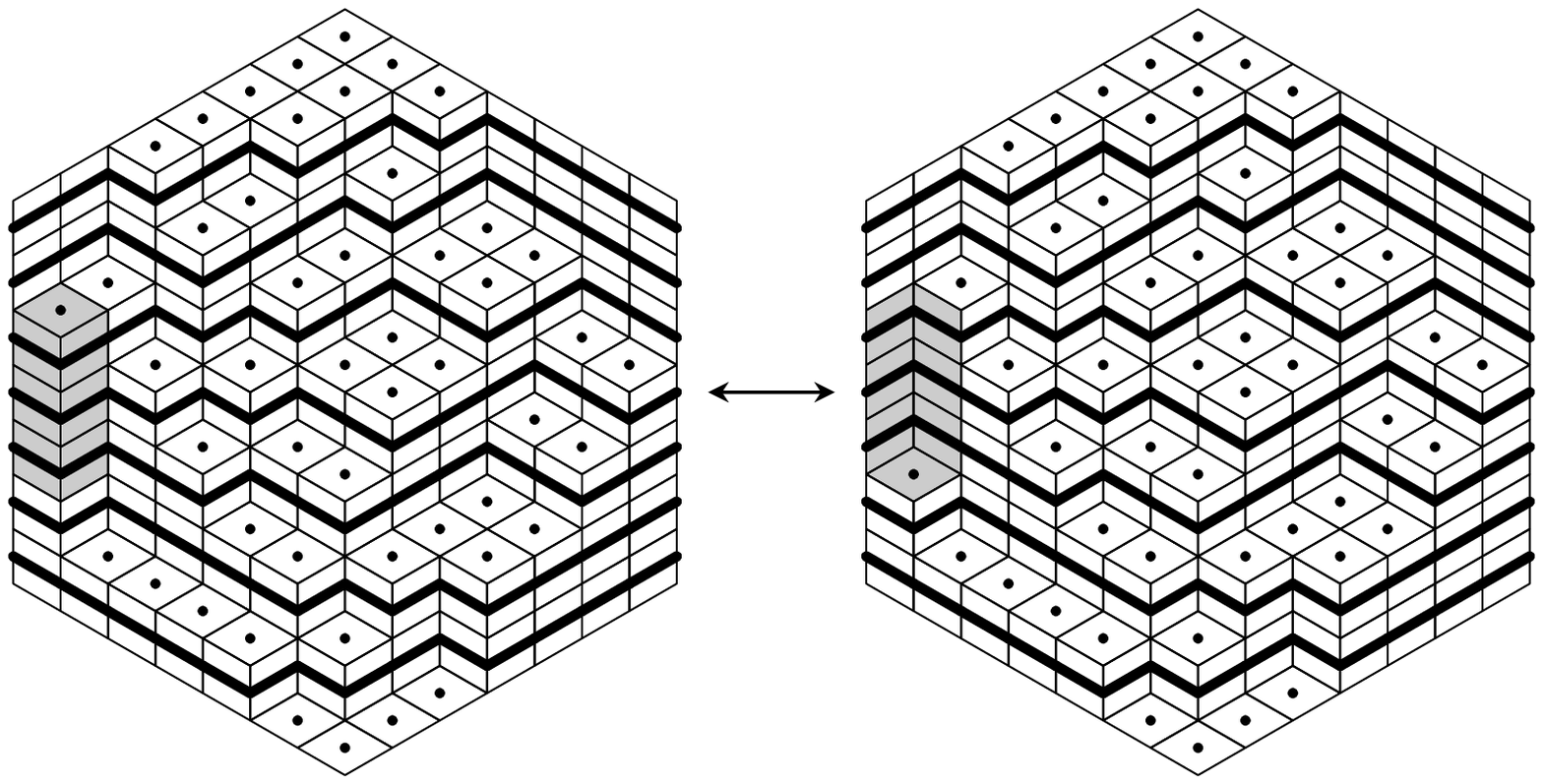}}
\caption{A tower move for the Luby-Randall-Sinclair Markov chain for lozenge tilings.  If the chain attempts to push up the third path from the bottom in the gray area, then it and the two paths blocking it are all three pushed up together with probability $1/3$.
}
\label{fig:tower}
\end{figure}
ergodic and symmetric (the probability that these $k+1$ paths get
pushed back to their original positions equals the probability that
they were moved in the first place), so the uniform distribution
remains the unique
stationary distribution.  And if there were no border effects, we
would be able to compute $E[\df]$, since when a site in a lattice path
is selected by the Markov chain, the expected change in the total
height function is determined by that point and its immediate
neighbors on the same path.  For hexagonal regions the borders cannot
obstruct the movements of the lattice paths, and we are able to get
both upper and lower bounds for these chains.  For general regions we
still obtain a good upper bound bound on the mixing time despite these
``border effects''.

The constraints that the borders of the region impose are that certain
locations of certain lattice paths have maximal or minimal values.
For instance, at the start and end of a lattice path, the maximal and
minimal values are identical.

\old{
Rather than work with the displacement function, we need to work with
the gap between an upper and lower lozenge tiling.  The advantage of
working directly with the gap function is that if the region's borders
affect the evolution of the gap, they can only make it on average
shrink faster.
}

Let $w$ be the width of the lozenge-tiling region in the lattice path
representation.  That is, $w$ is the distance between the leftmost
start of a path and the rightmost end of a path.  (For the path Markov
chain, $w$ was $n$.)  As in the single
path Markov chain, we find it convenient to center the region about
the origin, so that the $x$-co\"ordinates of the points on the paths
range from $-w/2$ to $w/2$.  For a given set of nonintersecting
lattice paths (a.k.a.\ 
routing), let $h_i(x)$ be the height of the $i$th path in the
routing at the given $x$-co\"ordinate.  If the path includes $x$ and
$x-1$, we have $h_i(x)=h_i(x-1)\pm 1/2$, so that pushing the path up
or down changes its height by 1.  For convenience, extend the
definition of $h_i$ to locations $x$ where the $i$th path does not
have a point, say by letting $h_i(x)$ take on its maximum possible
value consistent with the
constraint $h_i(x)=h_i(x-1)\pm 1/2$.  The ``displacement'' function of
$h$ that we will use is
\begin{equation}
\label{eqn:phi-lozenge}
 \Phi(h) = \sum_i \sum_{x=-w/2}^{w/2} h_i(x) \cos \frac{\beta x }{w}\ ,
\end{equation}
where $0\leq\beta\leq\pi$.  This function $\Phi$ is the natural generalization
of our earlier lattice path displacement function, but in the context of
lozenge tilings we put ``displacement'' in quotes because it does not measure
displacement from anything in particular.

Suppose that the Markov chain picks a particular location $x$ on the
$i$th path to try randomly pushing up or down.  Let $h'_i(x)$ be the
updated value of the gap at that location.  If the $i$th path does not
have an extremum at $x$, then $h'_i(x) = E[h'_i(x)] = h_i(x) =
[h_i(x-1)+h_i(x+1)]/2$, so 
\begin{equation}\label{eqn:edf}
E[\df(h)] = \left[\frac{h_i(x-1)+h_i(x+1)}{2} - h_i(x)\right] \cos
\frac{\beta x }{w}\ .
\end{equation}
Suppose instead that the $i$th path does have an extremum at $x$.
In the absence of border constraints or interactions
with other paths, the extremum is pushed up or down with probability
$1/2$, so that we would have $E[h'_i(x)] = [h_i(x-1)+h_i(x+1)]/2$.
Next we take into account interactions amongst the paths, but not
border effects yet.  Suppose the
Markov chain attempts to push the given local extremum in the opposite
direction that it is pointing, but that there are $k$ paths blocking
its movement.  Because we are using nonlocal moves,
with probability $1/(k+1)$ each of these $k+1$ paths is moved at
location $x$, with each affecting $\Phi(h)$ by $\cos(\beta x/w)$.
Thus Equation~\eqref{eqn:edf} continues to hold true. 

The only reason that \eqref{eqn:edf} might fail is if path $i$ has a
local maximum at $x$ and pushing it down violates the border
constraints (in which case $=$ in \eqref{eqn:edf} becomes $>$), or
the path has a local minimum at $x$ and pushing it up violates the
border constraints (where $=$ in \eqref{eqn:edf} becomes $<$).  For
some regions, such as the hexagon, all
the paths start at $-w/2$ and end at $w/2$, and the only border
effects are that the endpoints of the paths stay fixed.  For such
regions Equation~\eqref{eqn:edf} always holds for each path and each
$x$ such that $-w/2<x<w/2$, and it holds for $x=\pm w/2$ as well when
$\beta=\pi$.  Using this equality we will derive a lower bound on the
mixing time when the region is a hexagon.  For general regions it is
difficult to obtain a lower bound on the mixing time, but we still
obtain an upper bound.

\subsection{Mixing time upper bound}
\label{sec:lozenge-upper}

To get the upper bound we will work with a pair of routings with
heights $\hu_i$ and $\hl_i$, such that each path in the first routing
lies above the corresponding path in the second routing.
How we extended the definition of $\hu_i$ and $\hl_i$ to locations $x$
where path $i$ does not exist was a bit arbitrary, but we did it the
same way in both routings so that at these locations
$\hu_i(x)-\hl_i(x)=0$.
Then the gap function between the two routings is
$g=\hu-\hl$, and the gap is $\Phi(g) = \Phi(\hu)-\Phi(\hl)$, which is
zero when the routings are the same, and positive otherwise.

When we pick the same location $x$ on the same path in both routings,
and randomly push in the same direction, then from
Equation~\eqref{eqn:edf} we get
\begin{equation}\label{eqn:edg}
E[\df(g)] = E[\df(\hu)-\df(\hl)] = \left[\frac{g_i(x-1)+g_i(x+1)}{2} - g_i(x)
  \right] \cos \frac{\beta x }{w}
\end{equation}
unless the borders influence $E[\df(g)]$.  Suppose that the
borders influence $E[\df(\hu)]$ to be larger than the value given by
\eref{edf}.  Then $\hu_i(x)$ takes on its minimal possible value, and is
a local maximum, so $\hu_i(x\pm1)$ also take on their minimum possible
values.  But $\hu_i$ dominates $\hl_i$, so $\hl_i(x)$ and $\hl_i(x\pm 1)$ also
assume these same minimum possible values.  Thus the right-hand side
of \eqref{eqn:edg} is zero.  Since $\hu_i(x)$ and $\hl_i(x)$ are immobile, the
left-hand side is zero as well, so \eqref{eqn:edg} continues to be true.
Similarly, if the borders influence $E[\df(\hl)]$ to be smaller than the
value given by \eref{edf}, \eref{edg} continues to hold true.  If on
the other hand, the borders influence $E[\df(\hu)]$ to be smaller than
the value specified by \eref{edf}, and/or influence $E[\df(\hl)]$ to be
larger than that specified by \eqref{eqn:edf}, then we may replace the
second equality in \eref{edg} with $\leq$ to obtain a true statement.

Let $p$ denote the number of internal points on the paths in the
routing, i.e.\ the number of places where the Markov chain might try
to push a path up or down.  (The $p$ in this section was $n-1$ in the
section on the path Markov chain.)  Then when we stop conditioning on
a particular site of a certain path getting pushed, we use the same
derivation used in the proof of \lref{path-df} to conclude that $$
E[\df] \leq \frac{-1 + \cos(\beta/w)}{p} \Phi\ ,$$ with equality when
$\beta=\pi$, all paths start at $-w/2$ and end at $w/2$, and the only
restrictions on the locations of the paths are that that their
endpoints are pinned down and they do not intersect.  Then we use the
same argument used in the proof of \tref{path-upper} to find that the
Markov chain is within $\varepsilon$ of uniformity after $$
\frac{2+o(1)}{\beta^2} p w^2 \log \frac{m}{\fm\varepsilon}\ ,$$ where $m$ is
the number of (local) moves separating the upper and lower
configurations, and $\fm>\cos(\beta/2)\approx(\pi-\beta)/2$.  Taking
$\beta=\pi-\Theta(1/\log n)$ as before yields
\begin{theorem}
\label{thm:lozenge-upper}
For the Luby-Randall-Sinclair lozenge-tiling Markov chain on a region which has width $w$, has $m$ (local) moves separating the top and bottom configurations, and contains $p$ places where lattice path may be
moved, after
$$\frac{2+o(1)}{\pi^2} p w^2 \log \frac{m}{\varepsilon} $$
steps, coalescence will occur except with probability $\varepsilon$.
(For regions which are not simply connected, we mean coalescence within
a given connected component of the state space.)
\end{theorem}
For the bound stated in the introduction we used $p\leq n$ and
$m\leq 2 n^{3/2}$ \citep*{luby-randall-sinclair:markov-lattice}.
% in proof of thm 12 of journal version

For the hexagon of order $\ell$, $p=2\ell(\ell-1)$, $w=2\ell$, and $m=\ell^3$, so our
mixing time bound is $$\frac{48+o(1)}{\pi^2} \ell^4 \log \ell \ .$$
 % when $\varepsilon$ does not go to zero too quickly.

\subsection{Lower bound for the hexagon}
\begin{theorem}
For the regular hexagon with side length $\ell$, the mixing time of the
Luby-Randall-Sinclair Markov chain is at least $(8/\pi^2-o(1))\ell^4\log\ell$.
\end{theorem}
\begin{proof}
We apply \lref{anticonverge} here to lower bound the mixing time of
the lozenge tiling Markov chain proposed by \lrs, when the region is a
regular hexagon with side lengths $\ell$.  Our potential function has
the required contraction property, with $\gap\approx \beta^2/(2 p
w^2) \approx \pi^2 / (16 \ell^4)$.  $\fx = \ell^3/2$.  Since we are
using nonlocal moves, $\df$ can be as large as $\ell$.  Suppose that
when the Markov chain picks a site on a path and tries to push it,
there are $k$ paths in the way.  With probability $1/(k+1)$\ \ 
$\df=(k+1)\cos(\alpha)$, and otherwise $\df=0$.  Conditioning on this
site being selected, $E[(\df)^2]\leq k+1$, and in general $E[(\df)^2]
\leq R = \ell$.  Applying \lref{anticonverge}, we obtain a mixing time lower bound of
\begin{align*}
\frac{(16+o(1)) \ell^4}{\pi^2} \left[ 3 \log \ell + \frac{1}{2} \log
(\varepsilon/\ell^5)\right] &=\frac{8+o(1)}{\pi^2} \ell^4 \log \ell\ .\qedhere
\end{align*}
\end{proof}

\subsection{The other Luby-Randall-Sinclair Markov chains}
\label{sec:other-tower}

In addition to the Markov chain for lozenge tilings, \lrs introduced
Markov chains for domino tilings and Eulerian orientations.  
\citet{randall:98} has pointed out that our analysis in \sref{lozenge-upper}
is readily adapted more or less unchanged to their domino tiling
Markov chains.  It is much less obvious how to adapt the analysis to
the Eulerian orientation chains.  The reason for this difference is
that for the Eulerian orientation Markov chains, the nonlocal ``tower
moves'' overlap one another in a criss-cross fashion, whereas in the
lozenge tiling and domino tiling chains the towers are parallel to one
another.  We effectively gave each local move within a given tower the
same weight, and if we do the same for the Eulerian orientation chain,
we get the trivial weighting, which does not have the desired contraction
property.

\subsection{The local move Markov chains}

For a ``normal'' $\ell\times\ell$ region one might expect that for
typical configurations, the towers in the nonlocal tower moves are
fairly short, which suggests that while the Luby-Randall-Sinclair
Markov chain is much nicer to analyze rigorously, it may not be much
faster on these regions than the local moves Markov chain.  (However
for certain contrived regions, such as a pencil-shaped region
consisting of one long tower, the LRS Markov chain will be much faster
than the local-moves Markov chain.)
Thus we have a heuristic prediction that the local moves chain takes
$\Theta(\ell^4\log\ell)$ time to mix in ``normal'' $\ell\times\ell$ regions.
\citet{cohn:personal} tested this prediction by doing coupling time
experiments, and reported that it was ``about right.''  Then
\cite{henley:height} did some detailed heuristic calculations and
predicted that the relaxation time (reciprocal of the
spectral gap) was $\Theta(\ell^4)$ for a variety of models
that have associated height functions.
In an interesting recent development, \cite{destainville:rhombus}
experimented with the local moves chain for rhombus tilings
of octagonal regions where there are $6=\binom{4}{2}$ types of rhombuses,
and concluded that the $\Theta(\ell^4\log\ell)$
estimate holds for these tilings as well.

From a rigorous standpoint, \citet*{randall-tetali:local} established
a polynomial time upper bound on the mixing time of the local moves
chain.  Their approach was to use the mixing time bound from
\tref{lozenge-upper} on the nonlocal moves Markov chain to obtain a
bound on the chain's spectral gap, use techniques developed by
\citet*{diaconis-saloff-coste:compare-reversible} to compare the
spectral gaps of the local and nonlocal Markov chains, and then derive
a mixing time bound for the local moves chain from its spectral gap.
Their local moves mixing time bound was $O(n^2 w^2 h^2\log n)$, where
$n$ is the area, $w$ is the width, and $h$ is the height, or in other
words $O(\ell^8\log\ell)$ for $\ell\times\ell$ regions.
If rather than starting from our mixing time bound on the nonlocal
Markov chain, one instead starts from our bound on the spectral gap
(which was not explicitly given in an earlier version of this
article), then the $\log$-factors disappear from the mixing time bound
of the local moves Markov chain.

\section{The \kdk Markov chain}\label{sec:k&k}

As mentioned in the introduction, random generation of linear
extensions can be used to approximately count the number of linear
extensions of a partially ordered set, which is a \#P-complete problem
\citep*{brightwell-winkler:extension}.
\citet*{dyer-frieze-kannan:volume} showed one can generate an
(approximately) random linear extension of a partial order in
polynomial time, using a Markov chain on a certain polytope.
\citet*{matthews:extension} gave a different geometric Markov chain
for random linear extensions that runs in time $O(550 n^8 \log^3 n
\log 1/\varepsilon)$. % , and bounded the required precision.
\citet*{karzanov-khachiyan:extension} gave a combinatorial Markov
chain for linear extensions, and showed that it mixes in time $8 n^5
\log(|\Omega|/\varepsilon) \leq O(n^6 \log n)$, where $|\Omega|$ is the
number of linear extensions.  \citet*{dyer-frieze:volume} improved
the mixing time bound to $O(n^4 \log(|\Omega|/\varepsilon)) \leq O(n^5
\log n)$.  \citet*{felsner-wernisch:linear-extension} showed that the
\kdk Markov chain mixes in time $O(n^3 \log n)$ for a certain class
of partial orders, and that one can obtain an unbiased sample in this
time.  \citet*{bubley-dyer:extension} showed that a related Markov
chain mixes in time $O(n^3 \log n)$, and that the original \kdk Markov
chain mixes in time $O(n^3 \log n \log(|\Omega|/\varepsilon)) \leq O(n^4
\log^2 n)$.  We show here that the \kdk Markov chain mixes in time
$O(n^3 \log n)$, and exhibit a partial order for which the \kdk
Markov chain and \bd's variation of it both need order $n^3 \log n$
steps before they begin to get close to being random.

(Despite this progress, it still takes $O(n^5 \log^2 n\ 
\varepsilon^{-2} \log(n/\varepsilon))$ time to approximately count
linear extensions to within a factor $1+\varepsilon$
\citep*{bubley-dyer:extension}.  One can count the linear
extensions of series-parallel posets much more quickly, so the
aforementioned data-mining application
\citep*{mannila-meek:partial-order} restricted its attention to
these posets.)

\citet*{bubley-dyer:extension} use their simple yet powerful method
of path coupling \citep{bubley-dyer:p-couple} to bound the
mixing time of Markov chains related to the \kdk Markov chain for
random linear extensions.  In their generalization, the items at
positions $i$ and $i+1$ are considered with probability $f(i)$,
and the Markov chain transposes these items with probability
$1/2$, provided that doing so does not violate the partial order.  For
the \kdk \mc, $f(i)=1/(n-1)$ for $i=1,\ldots,n-1$.  Bubley and Dyer
show that if $f(i)$ is given by a parabola, $f(i)\propto i(n-i)$, then
the Markov chain mixes in time $(1/3+o(1))n^3\log n$, and then argue
using eigenvalue comparison techniques that the original \kdk chain
mixes in time no larger than $O(n^4\log^2 n)$.

We show here how to generalize \bd's analysis of these Markov chains,
and obtain an upper bound of $(4/\pi^2+o(1))n^3\log n$ for the \kdk
\mc.  We remark that $4/\pi^2$ is about $22\%$ larger than $1/3$, but
selecting a uniformly random location is easier than selecting
one according to a parabolic distribution.  We will see in section
\sref{sweep} that by doing updates in ``sweeps'' rather than at
independent uniformly random locations, the required number
transpositions can be cut in half.  So this analysis
marginally improves but does not significantly impact the time needed
to generate random linear extensions.  Mainly it serves to illustrate
the utility of the technique used throughout this article for analyzing
the mixing time of a variety of Markov chains that have been studied
before.

What we do is simply add weights to the distance function between linear
extensions.  If positions $i$ and $j>i$ are
transposed, \bd defined the ``width'' of the transposition to be $j-i$.
We will define the width to be $w(i,j)=\sum_{i\leq k<j}w(k)$, where the
$w(k)\geq 0$ are to be chosen later.  Given two linear extensions $X$ and
$Y$ of the partial order, a transposition sequence was defined to be a
sequence of linear extensions $X=Z_0,Z_1,\ldots,Z_r=Y$ such that $Z_k$ and
$Z_{k+1}$ differ by a single transposition.  The weight of a
transposition sequence is the sum of the widths of the transpositions,
and \bd define the distance $\delta(X,Y)$ between linear extensions $X$
and $Y$ to to be the weight of the minimum weight transposition
sequence.

\bd show in an appendix that (when each the $w(i)=1$) $\delta(X,Y) =
\delta_S(X,Y)$, where $\delta_S$ is Spearman's footrule, which is defined
by $\delta_S(X,Y)=(1/2)\sum_{i=1}^n |X(i)-Y(i)|$.  The proof is also
valid in the weighted scenario, when the definition of $\delta_S$ is
generalized to $\delta_S(X,Y) = (1/2)\sum_{i=1}^n
w(\min(X(i),Y(i)),\max(X(i),y(i)))$.

Given two permutations $A$ and $B$ that differ by single transposition
$(i,j)$ and which are updated to $A'$ and $B'$ using a coupled Markov
chain, \bd prove that
\begin{align*}
E[\delta(A',B')] &\leq \delta(A,B) + \frac{f(i-1)w(i-1)
  -f(i)w(i)-f(j-1)w(j-1) + f(j)w(j)}{2} \\
&= \delta(A,B) (1-\gap_{i,j})
\end{align*}
where
\begin{equation*}
\gap_{i,j} = \frac{1}{2} \frac{-f(i-1)w(i-1) +f(i)w(i) +f(j-1)w(j-1) - f(j)w(j)
}{w(i)+\cdots+w(j-1)}\ .
\end{equation*}
As always, they show this assuming constant weights $w$, but the same
proof holds for general positive weights $w$.

Letting $\gap=\min_{i,j}\gap_{i,j}$, their method of path
coupling gives an upper bound of $(1/\gap)\log(D/\varepsilon)$ on
the number of steps before the variation distance from stationarity is
at most $\varepsilon$, where $D$ is the ratio of the maximum distance
to the minimum positive distance.  Observe that $\gap_{i,j}=c$
(resp. $\gap_{i,j}\geq c$) for each $i$ and $j$ if and only if
$\gap_{i,i+1}=c$ (resp. $\gap_{i,i+1}\geq c$) for each $i$.

Given constant weights $w$, the choice of frequencies $f$ that
maximizes $\gap$ is given by a parabola.  Therefore \bd chose
$f(i)=i(n-i)/K$, where the normalizing constant is $K=(n^3-n)/6$.  It
is easily checked that $\gap_{i,j}=1/K$ when $j=i+1$, so this holds
for all $i$ and $j$.  With constant weights one can show $D\leq\lfloor
n^2/4\rfloor$, so their bound on the mixing time is $(1/3+o(1))n^3\log
n$.

Given constant frequencies $f$, the optimal choice of $w$ is
sinusoidal.  Let $w(i) = \cos(\beta(i/n-1/2))$ with
$0\leq\beta\leq\pi$; these weights are positive as required.  Since $$
\frac{-\cos(x-\delta) + 2\cos(x) - \cos(x+\delta) }{2 \cos(x)} =
1-\cos(\delta)\ ,$$ we have $$\gap_{i,j} \geq \frac{1-\cos(\beta/n)}{  n-1} \geq \frac{\beta^2}{2 n^3}$$ for
$j=i+1$ (we do not have equality when $i=1$ or $i=n-1$ since
$f(0)=f(n)=0$), so this bound holds for all $i$ and $j$.  
As in the proof of \tref{path-upper}, we take
%### $\beta=\pi-\Theta(\log\log n/\log n)$ 
$\beta=\pi-\Theta(1/\log n)$ so that $\gap$ is large, while
not making the ratio $D$ of the maximum distance to the minimum
positive distance too large.  Then the upper bound on the mixing time
is $$\frac{2 n^3}{\beta^2} \log \frac{D}{\varepsilon} = \frac{4+o(1)}{\pi^2}
n^3 \log n\ .$$

Of course when designing a Markov chain, we are free to pick both $f$
and $w$.  Optimizing them together would be an interesting challenge.

\old{
\citet*{bubley-dyer:extension} showed using their method
of path-coupling \citep*{bubley-dyer:p-couple} that this same
upper bound of $O(n^3 \log n)$ applies to the Markov chain applied to
any partially ordered set.  ... They were however unaware that there is a
matching lower bound (when the Markov chain is applied to the
antichain), and refered instead to the best previously published lower
bound of $\Omega(n^3)$ to argue that their analysis is reasonably tight.}

\section{Sweeps versus independent updates}\label{sec:sweep}

So far we have focused on updates where a random site is selected, and
then a local randomizing operation is performed at that site.  Often
in practice the various sites are updated in systematic
``sweeps'' rather than at
random.  For instance, for permutations or linear extensions, rather
than randomize a random adjacent pair of items, one may instead
randomize the items in positions $(1,2)\ (3,4)\ (5,6)\ \ldots$, and
then do positions $(2,3)\ (4,5)\ (6,7)\ \ldots$.  Likewise for lozenge
tilings, one may randomize the lattice paths at all places where the
$x$-co\"ordinate is even, then afterwards at all places where the
$x$-co\"ordinate is odd.  Call the first set of updates an even sweep,
and the second set an odd sweep.  In all cases where we have derived
upper bounds on the mixing time of a Markov chain (i.e.\
grand-coupling time for a random path or permutation by random
adjacent transpositions, grand-coupling time for
\lrs's chain on lozenge tilings, and the mixing time for the
\kdk Markov chain for random linear extensions) the same analysis that
worked for independent updates at each step also works for sweeps.
If one randomly chooses between even sweeps and odd sweeps, then we
have the same contraction property, with $\gap$ scaled up by a
factor of $(n-1)/2$ (or $p/2$ in the case of lozenge tilings).  The
mixing time bounds are then roughly the same though slightly better
than the bounds we would get when using the same total number of
transpositions (or path pushes) but at independent uniformly random
locations.  Successive even sweeps are redundant, as are successive
odd sweeps, so when we alternate, we perform about half as many moves (in continuous time)
to get the same value of the total variation distance or probability
of not coupling.

Note that we have not proved that the mixing time is actually twice
as fast, merely that our upper bound on it is half as large.

Perhaps more important than this factor of two savings is that many
fewer random bits are required to do the updates, since the locations
of the updates are deterministic.  For Markov chains with simple moves
such as these, generating pseudorandom bits can take an appreciable
fraction of the total running time.  Whether for this reason or for
simplicity, in practice the algorithms used to generate random tilings
have typically used systematic sweeps.

\section{Lattice paths and permutations revisited}
\label{sec:path-perm}

We have already given a quick-and-easy analysis of the adjacent-transposition Markov chain on lattice paths and on permutations, obtaining upper and lower bounds on the mixing time and coupling time that match to within constants.  In this section we give a more refined analysis which improves these constants.

\subsection{Upper bounds}
\label{sec:path-perm-upper}

Consider the Markov chain on permutations which exchanges a random
adjacent pair with probability $1/2$, and the pairwise coupling for
which the choice of adjacent pair is always the same in the two Markov
chains, and the decision of whether or not to exchange is also the
same in both chains, unless an exchange in one chain but not the other
would decrease the Hamming distance between the two permutations, in
which case the exchange is done in a random one of the two
permutations but not the other.  
This coupling was also used by \cite{aldous:group-walks}.
Let us focus on how a given item
moves in the two permutations.  The state space is the $n\times n$
grid, representing the location of the item in the two different
permutations.  A typical state $(x,y)$ transits to its four neigboring
states $(x\pm 1,y)$ and $(x,y\pm 1)$, each with probability
$\alpha=1/(2(n-1))$, with the following exceptions: (1) if $x=y$ then
the transitions are to $(x+1,y+1)$ and $(x-1,y-1)$, each with
probability $\alpha$, and (2) if the transition would be to a pair
outside the $n\times n$ grid, the transition instead self-loops.
Observe that the $x$-coordinate is a simple random walk on a chain of
length $n$, and similarly for the $y$-coordinate.  Furthermore, any
state $(x,y)$ with $x\neq y$ is transient, so that eventually $x=y$.
This will take about $\Theta(n^2/\alpha)$ steps.  We will need a good
estimate of the probability that coalescence has not occured after a
large multiple of $\Theta(n^2/\alpha)$ steps, so we prove the following
lemma:
\begin{lemma}\label{lem:triangle}
After $T$ steps $\Pr[x_T\neq y_T] < 10 \exp[-T (1-\cos(\pi/n))/(n-1)]$.
\end{lemma}
Before proving this lemma, we derive from it our bounds on the coupling times.

\begin{theorem}
\label{thm:pairwise}
For the above pairwise coupling on permutations, the number of steps before the probability of coalescence is at least $1-\varepsilon$ is at most $(2/\pi^2+o(1))n^3 \log(10 n /\varepsilon)$.
\end{theorem}
\begin{proof}
After $T$ steps the probability that the two permutations are still different is at most the expected number of items that are in different positions in the two permutations, which is at most $10 n \exp[-T (1-\cos(\pi/n))/(n-1)]$.  Setting this equal to $\varepsilon$ gives the desired bound.
\end{proof}

For the permutation Markov chain, the $(2/\pi^2+o(1)) n^3 \log n$
upper bound for the variation distance threshold and the
$(4/\pi^2+o(1)) n^3 \log n$ upper bound for the separation distance
threshold follow immediately.  The same bounds for the lattice path
Markov chain follow from projecting the permutation to the lattice
path.

\begin{theorem}
\label{thm:path-perm-upper}
For the sort/reverse-sort grand coupling on permutations, the number of steps before the probability of coalescence is at least $1-\varepsilon$ is at most $(4/\pi^2+o(1))n^3 \log(10 n /\varepsilon)$.  For lattice paths the sort/reverse-sort grand coupling time is at most $(2/\pi^2+o(1))n^3 \log(10 n /\varepsilon)$.
\end{theorem}
\begin{proof}
What we have analyzed in \tref{pairwise} is a pairwise coupling of a
Markov chain on permutations, i.e.\ an update rule that updates pairs
of permutations.  This pairwise update rule does not extend to a grand
coupling, i.e., there is no update rule defined on all permutations
such that pairs of permutations evolve according to the above pairwise
coupling.  But let us look at the evolution of a threshold function
of the two permutations.  Normally both permutations are either sorted
or reverse-sorted at a location, which corresponds to a push-down or
push-up move in the threshold functions.  The exceptional case where
one permutation is sorted while the other is reverse-sorted occurs
when in one permutation items $i$ and $j$ are in adjacent locations $x$
and $x+1$, while the other permutation has items $j$ and $k$ at these
locations, and either $i<j<k$ or $k<j<i$.  Any given threshold
function will map $j$ to either $0$ or $1$, and then one of the two
permutations will have either an up-slope or down-slope at locations
$x$ and $x+1$, and for that permutation an observer would be unable to
tell whether a sort or reverse-sort operation was performed.  Thus
from the standpoint of an observer, it appears as if the two threshold
functions were evolving according to the monotone grand coupling
considered earlier, and our bound on the pairwise coupling time for
permutations translates to a bound on the grand coupling time for
lattice paths.

Next we can convert the bound on grand coupling time for lattice paths
into an upper bound on the grand coupling time for permutations for
the straightforward sort/reverse-sort coupling.  After $(4/\pi^2 +
o(1)) n^3 \log n$ steps the pairwise permutation coupling (and hence
the lattice path grand coupling) has coalesced except with probability
$\ll 1/n$, so the permutation grand coupling has coalesced except with
probability $\ll 1$.
\end{proof}

Surprisingly, experiments suggest rather
strongly that $4/\pi^2$ is in fact the correct constant, so that
essentially nothing was lost in converting the coupling time bounds
from permutations to lattice paths and back to permutations.

\begin{proof}[Proof of \lref{triangle}]
Since we are interested in the probability that the random walk has
not hit the diagonal, and the regions below and above the diagonal
behave symmetrically, let us consider the state transition matrix
$M_n$ for the random walk above the diagonal ($x<y$), where the random
walker reflects off the boundaries of the grid, and dies when it hits
the diagonal.  The matrix $M_n$ resembles a stochastic matrix, except
that for those rows corresponding to states next to the diagonal, the
row-sum will be less than one.  For the reader's convenience we
proceed to diagonalize the matrix $M_n$; other triangular regions
with different boundary conditions have been similarly
diagonalized in e.g.\ \citep*{kenyon-propp-wilson:temperley}.

For $0\leq j <k <n$, let the function $f_{j,k}(x,y)$ be defined by
$$f_{j,k}(x,y) = \cos(j\pi x/n) \cos(k\pi y/n) - \cos(j\pi y/n) \cos(k\pi x/n).
$$
For convenience let the values of the grid coordinates
 $x$ and $y$ range from $1/2$ to $n-1/2$.  Since
\begin{multline*}
f_{j,k}(x-1,y)+f_{j,k}(x+1,y)+f_{j,k}(x,y-1)+f_{j,k}(x,y+1)-4f_{j,k}(x,y)\\
  = [2\cos(j\pi/n)+2\cos(k\pi/n)-4]f_{j,k}(x,y)\ ,
\end{multline*}
$f_{j,k}$ is an eigenvector of the nearest-neighbor random walk on $\Z^2$
with transition probabilities $\alpha$, and its eigenvalue is
  $$\lambda_{j,k} = 1+\alpha[2\cos(j\pi/n)+2\cos(k\pi/n)-4]\ .$$
Since furthermore
$f_{j,k}(x,x)=0$,
$f_{j,k}(x,y)=f_{j,k}(-x,y)$, $f_{j,k}(x,y)=f_{j,k}(x,-y)$,
$f_{j,k}(x,y)=f_{j,k}(2n-x,y)$, and $f_{j,k}(x,y)=f_{j,k}(x,2n-y)$,
it follows that $f_{j,k}$ is also an eigenvector of $M_n$
with eigenvalue $\lambda_{j,k}$.

Next we show that any two of these $n(n-1)/2$ eigenvectors are
orthogonal, and that each is not identically zero:
\begin{align*}
\sum_{1/2\leq x < y \leq n-1/2} f_{j_1,k_1}(x,y) f_{j_2,k_2}(x,y)
  &= \sum_{x<y} \begin{aligned}[t]\big[
  &\cos(j_1\pi x/n) \cos(k_1\pi y/n) \cos(j_2\pi x/n) \cos(k_2\pi y/n) \\
 +&\cos(j_1\pi y/n) \cos(k_1\pi x/n) \cos(j_2\pi y/n) \cos(k_2\pi x/n) \\
 -&\cos(j_1\pi x/n) \cos(k_1\pi y/n) \cos(j_2\pi y/n) \cos(k_2\pi x/n) \\
 -&\cos(j_1\pi y/n) \cos(k_1\pi x/n) \cos(j_2\pi x/n) \cos(k_2\pi y/n)
\big]\end{aligned} \\
  &= \sum_{x,y} \begin{aligned}[t]\big[
  &\cos(j_1\pi x/n) \cos(k_1\pi y/n) \cos(j_2\pi x/n) \cos(k_2\pi y/n) \\
 -&\cos(j_1\pi x/n) \cos(k_1\pi y/n) \cos(j_2\pi y/n) \cos(k_2\pi x/n)
\big]\end{aligned} \\
  &=\begin{aligned}[t]
  &\sum_x\cos(j_1\pi x/n)\cos(j_2\pi x/n) \sum_y\cos(k_1\pi y/n)\cos(k_2\pi y/n) \\
 -&\sum_x\cos(j_1\pi x/n)\cos(k_2\pi x/n) \sum_y\cos(k_1\pi y/n)\cos(j_2\pi y/n)\end{aligned} \\
  &= \begin{aligned}[t]
    &\frac{1_{j_1=j_2}+1_{j_1=j_2=0}}{2} n \times
     \frac{1_{k_1=k_2}+1_{k_1=k_2=0}}{2} n \\
   -&\frac{1_{j_1=k_2}+1_{j_1=k_2=0}}{2} n \times
     \frac{1_{k_1=j_2}+1_{k_1=j_2=0}}{2} n\ .
  \end{aligned}
\end{align*}
Since $j_1<k_1$ and $j_2<k_2$ the second term is zero.
The first term is also zero unless both $j_1=j_2$ and $k_1=k_2$,
giving us orthogonality.  If $(j_1,k_1) = (j_2,k_2) = (j,k)$, then we find
that this inner product is $(1+1_{j=0}) n^2 / 4$,
so the eigenvectors are nontrivial.  Hence, we have an orthogonal eigenbasis
of the matrix $M_n$.

Suppose that the random walker starts at $(x_0,y_0)$.  Let
$\delta_{x_0,y_0}(x,y)$ be the function which is $1$ at the starting
location and zero elsewhere.  Let $J(x,y)$ denote the function which
is $1$ whenever $x<y$.  We have
\begin{align*}
\Pr[x_T\neq y_T]
 &= (\delta_{x_0,y_0} M_n^T) \cdot J \\
 &= \left(\sum_{j<k} \frac{\delta_{x_0,y_0}\cdot f_{j,k}}
                          {f_{j,k}\cdot f_{j,k}}       f_{j,k} M_n^T \right)\cdot
    \left(\sum_{j<k} \frac{J\cdot f_{j,k}}
                          {f_{j,k}\cdot f_{j,k}}       f_{j,k}\right)     \\
 &=       \sum_{j<k} \frac{(\delta_{x_0,y_0}\cdot f_{j,k}) (J\cdot f_{j,k})}
                          {f_{j,k}\cdot f_{j,k}}       \lambda_{j,k}^T    \\
 &<       \sum_{j<k} \frac{(2)\left(2\binom{n}{2}\right)}
                          {(1+1_{j=0})n^2/4}           \lambda_{j,k}^T    \\
 &<  8 \sum_{j<k} \lambda_{j,k}^T\ .
\end{align*}
We need a bound on $\lambda_{j,k}$ to bound this summation, and to this end
consider the line passing through $(0,\cos 0)$ and $(t,\cos t)$.  When
$t=\pi/2$, the line is at least as high as $\cos s$ when $t\leq s \leq
\pi$.  If $t<\pi/2$, the line's slope increases towards $0$, so it
continues to be above $\cos s$ when $\pi/2 \leq s \leq \pi$, and by
concavity of $\cos s$ for $0\leq s \leq \pi/2$, the line is also above
$\cos s$ when $t\leq s \leq \pi/2$.  Taking $t=\pi/n$ (assume $n\geq2$ ---
the lemma is trivial if $n=1$) and $s=j\pi/n$, we have
\begin{align*}
\cos(j\pi/n)    &\leq 1-\frac{j\pi/n}{\pi/n}(1-\cos(\pi/n)) \\
2\cos(j\pi/n)-2 &\leq -2j(1-\cos(\pi/n)) \\
2\cos(j\pi/n)+2\cos(k\pi/n)-4 &\leq -2(j+k)(1-\cos(\pi/n)) \\
\lambda_{j,k}   &\leq 1-\alpha 2(j+k)(1-\cos(\pi/n))\ .
\intertext{Let $c=2\alpha(1-\cos(\pi/n)) \approx \alpha\pi^2/n^2$ so that}
\lambda_{j,k}   &\leq 1-(j+k)c < \exp[-(j+k)c]\ .
\intertext{Since $\alpha=1/(2(n-1))$, $\lambda_{j,k}\geq 1-8/[2(n-1)]\geq 0$ when $n\geq 5$,
and $\lambda_{j,k}\geq 0$ for $n=2,3,4$ as well,
so we may take the $T$th power of both sides to get}
\lambda_{j,k}^T &\leq \exp[-(j+k)c T] \\
\sum_{k=j+1}^{\infty}\lambda_{j,k}^T
                &\leq \frac{\exp[-(2j+1)c T]} {1-\exp(-c T)} \\
\Pr[x_T\neq y_T] <
8 \sum_{0\leq j<k} \lambda_{j,k}^T
                &\leq 8 \frac{\exp[-c T]} {[1-\exp(-2 c T)][1-\exp(-c T)]\ .}
\intertext{The lemma is trivial unless $\exp[-c T]\leq 1/10$, in which case
$ 8/[1-\exp(-2 c T)]/[1-\exp(-c T)] \leq 8\cdot 100/99 \cdot 10/9 = 8000/891<10$, so that}
\Pr[x_T\neq y_T] &< 10 \exp[-c T]\ . \qedhere
\end{align*}
\end{proof}

\subsection{Lower bounds}
\label{sec:path-perm-lower}

\begin{theorem}
\label{thm:path-lower-c}
For the Markov chain on lattice paths in a $n/2 \times n/2$ box, the time
it takes the top path and bottom path to coalesce is with high probability
at least $(1-o(1))2/\pi^2 n^3 \log n$.
\end{theorem}

\begin{proof}
As an upper lattice path $\hu$ and lower lattice path $\hl$ evolve
together via the push down / push up coupling, let us look at the
difference path $h=\hu-\hl$.  If $\hu$ goes up and $\hl$ goes down,
which we will denote \pp{U}{D}, then the difference path $h$ goes up,
which we denote with \U.  If $\hu$ goes down and $\hl$ goes up
(\pp{D}{U}), then $h$ goes down (\D).  In the remaining two cases
(\pp{U}{U} and \pp{D}{D}) the difference path remains flat (\F).  We
may view the difference path as a string of \U, \F, and \D particles,
and it is easy to check that the evolution of the difference path is a
Markov process: If the particles at the updated site are
\p{UU}=\pp{UU}{DD}, then they remain \pp{UU}{DD}=\p{UU}.  If a
\p{UD}=\pp{UD}{DU} is updated, the result is either \pp{UD}{UD}=\p{FF}
or \pp{DU}{DU}=\p{FF}.  If a \p{UF} is updated, the underlying paths
might be \pp{UU}{DU} and then change to \pp{UU}{UD}=\p{FU} or
\pp{UU}{DU}=\p{UF}, or the underlying paths might be \pp{UD}{DD} and
then change to \pp{UD}{DD}=\p{UF} or \pp{DU}{DD}=\p{FU}.  Likewise, if a
\p{FF} is updated, there are four possibilities for the underlying paths,
and in each case the updated configuration is \p{FF}.  The other cases
(\p{DD}, \p{DU}, \p{DF}, \p{FU}, and \p{FD}) are similar and related
to the above cases by symmetry.  We may summarize the update rules for
the string of \U's, \D's, and \F's as follows: pick a random adjacent
pair, and with probability $1/2$ exchange them; when a \D and \U are
exchanged past each other, they both turn into \F's.  If we start with
the top path and bottom path, then in the difference path every \U
will be to the left of every \D.

We will do a number of comparisons between a random permutation $\sigma$ and
the difference lattice path $h$.  For the $k$th comparison ($0\leq k\leq
n$), look at the locations of cards $1,\ldots,k$, and in particular
their relative order.  Let $\tau_\U(1)$ denote the first of these cards
encountered in a left-to-right scan, and in general $\tau_\U(i)$, $1\leq
i\leq k$, denotes the $i$th such card encountered.  Label the first
$k$ U's of the difference path with the numbers
$\tau_\U(1),\ldots,\tau_\U(k)$.  Similarly let $\tau_\D(i)$, $1\leq i\leq k$,
denote the $i$th card from the cards $n+1-k,\ldots,n$ to be
encountered in a right-to-left scan of the permutation, and label the
$i$th to last \D of the difference path with $\tau_\D(i)$.  We leave the
remaining particles of the difference path unlabeled.  When we evolve
the difference path via random exchanges, we will let labeled
particles be exchanged past each other, but a labeled \U or \D may not
be exchanged past an unlabeled \U or \D.  This rule for the labels does
not affect the evolution of the unlabeled difference path, but it is
important for our understanding of it.

Initially for each $1\leq i\leq k$, the position of card $\tau_\U(i)$ in
the difference path is weakly to the left of card $\tau_\U(i)$ in the
permutation, while the position of card $\tau_\D(i)$ in the difference
path is weakly to the right of card $\tau_\D(i)$ in the permutation.
We will pick the same random adjacent pair in the permutation as in
the labeled difference path, and make the same decision as to whether
or not to exchange the adjacent items.  Consider the first time that
the above invariant fails to hold, say that card $\tau_\U(i)$ in the
labeled difference path moves to the right of card $\tau_\U(i)$ in the
permutation.  On the previous step, card $\tau_\U(i)$ was in the same
location in the difference path and the permutation.  The exchange
could not have been to the right of the card $\tau_\U(i)$, because
exchanges in the permutation always succeed, nor could the exchange be
to the left of card $\tau_\U(i)$, as any particle to the left of card
$\tau_\U(i)$ is either an \F or a labeled \U, and such exchanges succeed.
Thus the invariant is maintained.

Consider the locations of two cards $i$ and $j$ in a random permuation
$\sigma$, or two labels in the difference path $h$.  Let the weighted
gap between them be defined by $\wgap(i,j) = \sum_x \sin(\pi x/n)$,
where the sum is taken over positions $x$ between the two cards, and
is negative if card $j$ occurs before card $i$.
Within a random permutation $\sigma$ we have
\begin{align*}
 E[|\wgap_\sigma(i,j)|]
   &= \frac{\sum_{x=1}^{n-1} x(n-x)\sin\frac{\pi x}{n}}{\binom{n}{2}} \\
   &\approx 2 n \int_0^1 u(1-u)\sin(\pi u) du = \frac{8}{\pi^3} n\ .
\end{align*}

The area under the difference path is the sum of the locations of
the \D particles minus the sum of the locations of the \U particles.
The potential function (the weighted area) is
\begin{align*}
 \Phi(h) &=    \sum_{x=0}^n h(x) \sin\frac{\pi x}{n} \\
         &=    \sum_{i=1}^{\max_x h(x)} \wgap_h(\text{$i$th \U},\text{$i$th \D}) \\
         &\geq \sum_{i=1}^{k} \wgap_h(i,n+1-i)
\intertext{if $k\leq\max_x h(x)$.  As $\wgap_h(i,n+1-i)\geq \wgap_\sigma(i,n+1-i)$ and also
$\wgap_h(i,n+1-i)\geq 0$ whenever $i\leq \max_x h(x)$, we have}
 \Phi(h) &\geq \sum_{i=1}^{k} \max(0,\wgap_\sigma(i,n+1-i))
\intertext{for any $k\leq\max_x h(x)$, so}
 \Phi(h) &\geq \sum_{i=1}^{n/2} 1_{i\leq\max_x h(x)} \max(0,\wgap_\sigma(i,n+1-i))\ .
\end{align*}
Since we started with a random permutation and the dynamics are reversible,
then even conditional upon all the past moves, the permutation is still uniformly random.  In particular $\wgap_\sigma(i,n+1-i)$ is independent of the maximum height $\max_x h(x)$ of the difference path, so that
\begin{align*}
E[\Phi(h)] &\geq \sum_{i=1}^{n/2} \Pr[i\leq\max_x h(x)] E[\max(0,\wgap(i,n+1-i))] \\
E[\Phi(h)] &\geq E[\max_x h(x)] [(8/\pi^3+o(1)) n / 2] \\
E[\max_x h(x)] &\leq (1+o(1))\frac{\pi^3}{4 n} \Phi(h_0) (1-\gap)^t
\intertext{and since $\Phi(h_0)\doteq (2/\pi^2)n^2$}
E[\max_x h(x)] &\leq (1+o(1))\frac{\pi}{2} n (1-\gap)^t\ .
\end{align*}
Note that this gives another proof that coalescence is likely after $t =
(2/\pi^2+o(1)) n^3 \log n$ steps.

Notice that the difference path never changes by more than one at a
time, and only if a \U or \D particle moves.  There are $2\max_x h(x)$
\U and \D particles, each particle can move in one of two directions, and a
given proposed exchange occurs with probability $1/2$.  Thus $$
E[\Delta\Phi^2|h(\cdot)] \leq 2\max_x h(x)/(n-1)\ .$$

\begin{align*}
E[\Phi^2(h_t)|h_{t-1}]  &= (1-2\gap)\Phi^2(h_{t-1}) + E[\Delta\Phi^2|h_{t-1}(\cdot)]\\
E[\Phi^2(h_t)]  &\leq (1-2\gap)E[\Phi^2(h_{t-1})] + E[E[\Delta\Phi^2|h_{t-1}(\cdot)]]\\
E[\Phi^2(h_t)]  &\leq (1-2\gap)E[\Phi^2(h_{t-1})] + (\pi+o(1)) (1-\gap)^{t-1}
\intertext{where here the $o(1)$ term depends only on $n$.  By induction}
E[\Phi^2(h_t)]  &\leq (1-2\gap)^t \Phi^2(h_0) + (\pi+o(1)) (1-\gap)^t/\gap \ .
\intertext{Subtracting $E[\Phi(h_t)]^2 = (1-\gap)^{2t} \Phi^2(h_0)$,}
\Var[\Phi(h_t)] &\leq (\pi+o(1)) (1-\gap)^t/\gap \\
\Var[\Phi(h_t)] &\leq (2/\pi+o(1)) n^3 (1-\gap)^t \\
\Var[\Phi(h_t)] &\leq (\pi+o(1)) n \Phi(h_t)\ .
\end{align*}

Thus if $\Phi_t \gg \pi n$, w.h.p.\ $\Phi(h_t)>0$, so that the time until
coalescence is w.h.p.\ at least $(1-o(1))2/\pi^2 n^3 \log n$.
\end{proof}

\section{Exclusion and interchange processes}\label{sec:exclusion}

In this section we show how to apply \lref{anticonverge} to lower bound
the convergence rate of exclusion and interchange processes.
Several of these Markov chains where studied by
\citet*{diaconis-saloff-coste:comparison}, who derived upper bounds on
their mixing times but did not have lower bounds that matched to
within constant factors.  The lower bounds derived here match (most of) their
upper bounds to within constant factors.

The interchange process describes particles moving around on an
undirected (but possibly weighted) graph.  At each time step, a random
edge of the graph is selected, and the particles at either endpoint of
that edge are exchanged.  The particles could be 0's and 1's, or they
could, for instance, be distinct numbers from $1$ up to the number $n$
of vertices.  
\old{If there are $k$ particles of one type, and $n-k$
particles of another type, then the mixing time lower bound for the
case of $k$ 1's and $n-k$ 0's is also a lower bound on the mixing time
for the given set of particles, so we can restrict are attention to
the case where the particles are 0's and 1's.}

\cite{lee-yau:log-sobolev} studied the logarithmic-Sobolev constants
and the $L_2$-mixing times of some exclusion and exchange processes.
The exclusion process may also be viewed as the infinite-temperature
limit of Kawasaki dynamics for the Ising model, see e.g.\ 
\citep{cancrini-martinelli:kawasaki} or \citep{lu-yau:kawasaki}.

Let $Q_{x,y}$ denote the probability that the Markov chain exchanges
the particles at locations $x$ and $y$, and for convenience let
$Q_{x,x}=1-\sum_{y\neq x} Q_{x,y}$.  Then $Q_{x,y}$ is the
state-transition matrix for the location of a particular particle.
Let $v$ be a right-eigenvector of matrix $Q$ (the matrix is symmetric,
so the left and right eigenvectors are the same) with eigenvalue
$1-\gap$.  The lemma requires $\gap>0$, but in general one would
expect that eigenvectors with smaller $\gap$'s will give better
lower bounds.

For convenience let us assume for the moment that all the particles are
distinguishable, so that the state of the Markov chain is a
permutation $\sigma$ of size $n$.  Define $$\Phi(\sigma) = \sum_{i=1}^k
v_{\sigma(i)}\ ,$$ and $\fx = \max_\sigma \Phi(\sigma)$.
After one step of the Markov chain we have
\begin{align*}
E[\Phi(\sigma')|\sigma]
  &= \sum_{i=1}^k E[v_{\sigma'(i)}|\sigma] \\
  &= \sum_{i=1}^k \sum_y Q_{\sigma(i),y} v_y\\
  &= \sum_{i=1}^k (1-\gap) v_{\sigma(i)}\\
  &= (1-\gap) \Phi(\sigma)
\end{align*}
so that this $\Phi$ has the contraction property required by
\lref{anticonverge} (assuming $0<\gap\leq 2-\sqrt{2}$).
In effect we used an eigenvector of the graph to define an eigenvector
of the interchange process.  For further information comparing the
eigenvalues of the graph with those of the interchange process, see
e.g.\ \cite{handjani-jungreis:interchange}.

With $\df = \Phi(\sigma')-\Phi(\sigma)$ denoting the change in $\Phi$ that
occurs after one step of the Markov chain, let $R$ be an upper bound
on the largest value that $E[(\df)^2]$ can take.  In general we can
take $$R\leq \max_{\text{$x,y$: $Q_{x,y}>0$}} (v_x-v_y)^2\ ,$$ but in
some cases we might find a better bound.

In general it need not be the case that all the particles are
distinguishable.  If there is a set $A$ of $k$ particles
that are distinguishable from the remaining $n-k$ particles, then we
can still define $$ \Phi(\text{state}) = \sum_{x: \text{particle of type $A$ 
at location $x$}} v_x\ .$$  We arbitrarily label the $A$ particles
$1,\ldots,k$ and the remaining particles $k+1,\ldots,n$.  The
evolution of $\Phi$ is exactly the same as it was when all the
particles were distinguishable.

One may add self-loops to the Markov chain to avoid periodicity
problems --- say the probability of a nontrivial transition is
$\alpha$.  Typical choices for $\alpha$ are $\alpha=1/2$ or
$\alpha\rightarrow 0$ (continuous time).

\subsection{Shuffling cards on a hypercube}

Here we lower bound the mixing time of the Markov chain considered by \dsc
that shuffles
cards via random transpositions where each transposition is an edge of
the hypercube.  The underlying graph from which we need an eigenvector
is the hypercube $\Z_2^d$, and the state space of the Markov chain is
$\Z_2^d !$ or $\binom{\Z_2^d}{2^{d-1}}$ depending on whether we want
to shuffle distinct particles or $2^{d-1}$ 0's and $2^{d-1}$ 1's (the
same lower bound applies to both cases).
(Here $V!$ denotes the set of permutations of a set $V$, and
$\binom{V}{n}$ denotes the set of subsets of $V$ containing $n$ items.)

We can take our eigenvector to be the function that is $1$ on those
vertices of the hypercube whose first co\"ordinate is 0, and $-1$ on
the other vertices.  When we follow a particular particle, the
probability that the particle gets moved across the first coordinate
is $\alpha/(d 2^{d-1})$, so our function is an eigenvector for which
$\gap = \alpha/(d 2^{d-2})$.  Here $\fx = 2^{d-1}$, and our bound
on $R$ is $4\alpha$.  Substituting into \lref{anticonverge}, we obtain a
lower bound of $$ (1-o(1)) \frac{d 2^{d-2}}{\alpha} \left[\log 2^{d-1} +
  \frac{1}{2} \log \frac{\varepsilon }{16 d 2^{d-2}} \right] =
(1-o(1)) \frac{\log 2 }{8 \alpha} d^2 2^d$$
for bounded values of $\varepsilon$.
It appears that this bound is correct up to constant factors.
%%
%%  The upper bound obtained by \dsc was order $d^2 2^d$,
%%  so this bound is correct up to constant factors.
%%
%  Saloff-Coste reports that the stated upper bound should have been $d^3 2^d$

\subsection{High-dimensional product graphs}

For a fixed connected graph $G$, consider the nearest neighbor random
walk on $G^d$.  We can imagine that there is a particle in each of $d$
disjoint copies of the graph $G$, where the particle in the $i$th copy
of $G$ gives the $i$th co\"ordinate of the walker.  At each time step a
particle chosen from a random copy of $G$ makes a move.  For example,
the Ehrenfest urn model from statistical mechanics is essentially
random walk on $\Z_2^d$, so here $G$ consists of two vertices and an
edge.  A common choice for the factor $\alpha$ by which to slow down
the walk is $\alpha=d/(d+1)$, in addition to the usual $\alpha=1/2$ and
$\alpha\rightarrow 0$.

The mixing threshold for random walk on $G^d$ has already been
determined, so it is an instructive exercise to check that
\lref{anticonverge} gives a sharp lower bound in this case.
For $\Z_2^d$ \citet*{diaconis-shahshahani:cube-convergence} showed
that there is a sharp variation mixing threshold at $(1/4\pm o(1)) d
\log d$ steps.  \citet*[Chapt.~7, sect.~1.7]{aldous-fill:book} state
that the same approach works for $G^d$ for more general graphs $G$.
In continuous time,
\citet*[Theorem~2.9]{diaconis-saloff-coste:log-sobolev} prove an upper
bound on the variation mixing time for random walk on $G^d$.  The
lower bound that we get from \lref{anticonverge} differs from the
upper bound by a factor of $1-o(1)$.

\citet*[sect.~7]{aldous-diaconis:walk-times} determined the mixing time
threshold of a related random on $G^d$ where at each time step, the
particles in each copy of $G$ get moved; we note that the lower bound
from \lref{anticonverge} is tight in this case as well.

To lower bound the mixing time of random walk on $G^d$, the underlying
graph for which we need an eigenvector is $\overbrace{G+\cdots+G}^d$.
Suppose that we have an eigenvector $v$ of $G$ with eigenvalue
$1-\gap_0$, where $\gap_0$ is the spectral gap.  We take our
eigenvector to be the canonical extension of $v$, i.e.\ the function
that assigns to each vertex $x$ of $G+\cdots+G$ the value of the
eigenvector $v$ in the copy of $G$ in which $x$ resides.  The
probability that a particular particle gets moved at a time step is
$\alpha/d$, so our function is an eigenvector for which $\gap =
\alpha\gap_0/d$.  (For the walks considered in \citep*{aldous-diaconis:walk-times} we have $\gap = \alpha\gap_0$.)
Here $\fx = \Theta(d)$, and our bound on $R$ is
$\Theta(\alpha)$.  (For the walks considered in \citep*{aldous-diaconis:walk-times} the bound on $R$ is $\Theta(\alpha d)$.)
Substituting into
\lref{anticonverge}, we obtain a lower bound of
\begin{equation}
\label{eqn:walk1}
\frac{ \log \Theta(d) + \frac{1}{2} \log
  \frac{\varepsilon\alpha\gap_0 }{d\Theta(\alpha)} }{\log 1/(1-\alpha\gap_0/d)} = (1-o(1)) \frac{d \log d }{2\alpha\gap_0}
\end{equation}
when $\varepsilon$ does not go to zero too quickly.

\remark
This application of \lref{anticonverge} is the natural generalization
of the approach that \citet*{diaconis-shahshahani:cube-convergence}
used to get their lower bound on the mixing time of $\Z_2^d$.  For the
walk on $\Z_2^d$, the eigenvalues of $\Z_2$ are $1$ and $-1$, so
$\gap_0 = 2$, which upon substition into lower bound \eqref{eqn:walk1}
gives the familiar $\frac14 d\log d$.  For $\Z_2^d$ the eigenvector
used above just counts the difference between the number of ones and
the number of zeros, which is the same test function that
\citet*{diaconis-shahshahani:cube-convergence} used.

\subsection{Shuffling cards on a grid} \label{sec:grid}

Here the cards (or 0's and 1's) are arranged on an $\ell\times m$ grid.
The state space is $([\ell]\times [m])!$ or
$\binom{[\ell]\times [m]}{ \lfloor \ell m/2 \rfloor}$, and the graph
for which we need an eigenvalue is the $\ell \times m$ grid.  (Here
$[n]$ denotes $\{1,2,\ldots,n\}$, and we identify the vertices of the
grid with $[\ell]\times[m]$.)  When $\ell=m$, \dsc show that order
$(\ell m)^2 \log (\ell m)$ steps suffice for the Markov chain to equilibrate,
and conjecture that
this is the correct order of magnitude (but see \sref{dsc}).
%  We shall confirm this conjecture.

Our asymptotic results will be valid as $\max(\ell,m)$ gets large, 
for convenience suppose that $\ell\geq m$.  Suppose the vertices are
labeled as $(i,j)$ for $1\leq i\leq \ell$ and $1\leq j\leq m$.  We can
take our eigenvector to be the function that is $\cos(\pi(i-1/2)/\ell)$
at vertex $(i,j)$.  There are $E=\ell(m-1)+m(\ell-1)$ edges of this grid
graph.  One can verify that this function is an
eigenvector with eigenvalue $$1-\gap = 1-2\alpha/E + (\alpha/E) 2
\cos(\pi/\ell)\ ,$$  so $\gap \sim \alpha\pi^2 / (E \ell^2)$.  The bound
for $R$ is order $\alpha(1/\ell)^2$, and $\fx$ is order 
$\ell m$.  Applying \lref{anticonverge}, we obtain a mixing time lower bound of
$$ (1-o(1)) \frac{E \ell^2 }{\alpha\pi^2} \left[ \log \Theta(\ell m) 
  + \frac{1}{2} \log \Theta(\varepsilon/E)\right]
 = (1-o(1)) \frac{\ell^2 (\ell-1/2) (m-1/2) }{\alpha\pi^2} \log (\ell m) $$
when $\varepsilon$ does not get too small too quickly
as $\ell m$ gets large.

\remark
We can recover our mixing time lower bound for shuffling cards by
adjacent transpositions by substituting $m=1$ and $\alpha=1/2$
into our lower bound for the $\ell\times m$ grid.

\subsection{\dsc's grid shuffling process} \label{sec:dsc}

The actual Markov chain that \dsc considered had transposition probabilities
that were slightly higher along edges that touched the border of the
grid.  This was because their update rule was to pick a random vertex,
and exchange the particles along a random edge incident to that
vertex.  Thus each edge $(u,v)$ is selected with probability
proportional to $1/d(u) + 1/d(v)$, where $d(u)$ and $d(v)$ are the
degrees of its endpoints.  In some sense it is clear that the slight
non-uniformity in the probability with which we select edges can not
affect the mixing time too much, but in the introduction we promised
a rigorous lower bound for \dsc's Markov chain, so we shall supply one.
Finding an explicit
eigenvector given these boundary conditions seems somewhat painful,
as does approximating one with sufficient accuracy, so
we take a different approach.  We show how to obtain mixing time lower
bounds using only an approximate contraction property.

Rather than try to approximate an eigenvector, it is easier to make
use of the fact that we have an exact eigenvector $\Phi$ to an
approximate state transition matrix, namely the state transition
matrix considered in \sref{grid} and the eigenvector $\Phi$ that
we used there
 $$\Phi(S) = \sum_{\text{$i,j$: state $S$ has particle at $(i,j)$}} \cos(\pi(i-1/2)/\ell)\ .$$
(For convenience we slow down the chain in \sref{grid} by a factor of
$\alpha=E/(2\ell m)=1-1/(2\ell)-1/(2m)$ when doing the comparisons,
so that the probability of a transition occuring on a given edge is
simply $1/(2\ell m)$.)
After one update of the approximate state transition matrix we have
$ E[\Phi(S')|S] = (1-\gap)\Phi(S)$
with $\gap$ as given above.
But with the actual state transition matrix, there is a
$O((\ell+m)/(\ell m))$ chance that a vertex on the boundary will be
selected, causing the exchange process to do something different than
the approximate exchange process.  To bound
$$ \left| E[\Phi(S')|S] - (1-\gap)\Phi(S)\right|$$
let us focus on a single particle at a time
and then make use linearity of expectations and the triangle inequality.
If the particle is near the boundary, the corresponding difference is at
most $O(1/m \ell^3)$; if the particle is not near the boundary then the
corresponding difference is $0$.  Since there are $O(\ell+m)$ particles near
the boundary,
with $\delta = O(1/\ell^3 + 1/(m\ell^2))$
 we have
$$ E[\Phi(S')|S] = (1-\gap)\Phi(S) \pm \delta\ .$$
By induction
$$E[\Phi(S_t)|S_0] = (1-\gap)^t \Phi(S_0) \pm \frac{\delta}{\gap}\ .$$
Likewise
\begin{align*}
E[\Phi^2(S_{t+1})|S_t]
 &= \Phi^2(S_t) + 2\Phi(S_t) E[\df|S_t] + E[(\df)^2|S_t]\\
 &= (1-2\gap) \Phi^2(S_t) \pm 2 \Phi(S_t) \delta + E[(\df)^2|S_t] \\
E[\Phi^2(S_{t+1})] 
 &= (1-2\gap) E[\Phi^2(S_t)] \pm 2 \Phi(S_0) \delta (1-\gap)^t \pm \frac{2\delta^2}{\gap} + E[(\df)^2]
\intertext{and so by induction}
E[\Phi^2(S_{t})] 
 &\leq (1-2\gap)^t\Phi(S_0) + \frac{2\Phi(S_0)\delta}{\gap}(1-\gap)^t +
       \frac{\max E[(\df)^2]}{2\gap} + \frac{\delta^2}{\gap^2}\ .
\intertext{Subtracting off $E[\Phi(S_t)]^2$ we get}
\Var[\Phi(S_t)]
 &\leq \frac{4\Phi(S_0)\delta}{\gap}(1-\gap)^t +
       \frac{\max E[(\df)^2]}{2\gap} + \frac{\delta^2}{\gamma^2}\ ,
\end{align*}
where in the last step we have assumed as in \lref{anticonverge} that
$1-2\gap$ is not excessively negative, i.e.\ that $0 < \gap \leq
2-\sqrt{2} \doteq 0.58$ (or else $0<\gap\leq 1$ and $t$ is odd).

For convenience let us follow \dsc in assuming $\ell=m$, so that
$\gap=(\pi^2/2+o(1))/\ell^4$ and $\delta=O(1/\ell^3)$.  We also have
$\Phi(S_0)=\Theta(\ell^2)$ when all the particles start on one side,
and $\max E[(\df)^2] = \Theta(1/\ell^2)$.
Let $K$ be a suitably large parameter to be selected in a moment.
Pick $t = \log(\ell/K) / \gap$, so that $E[\Phi(S_t)] = \Theta(K
\ell)$, and $\Var[\Phi(S_t)] = \Theta(K \ell^2)$.  But in stationarity
$E[\Phi(S)]=O(\ell)$ (indeed it is $0$) and $\Var[\Phi(S)] =
\Theta(\ell^2)$.  Thus for any given $\varepsilon$ we can take $K$
large enough so that, when we start the exclusion process with all the
particles on one side and run it for $\log(\ell/K(\varepsilon))/\gap$
steps, with probability $1-\varepsilon$ we are able to distinguish the
configuration from a state drawn from stationarity.  This gives us the
desired mixing time lower bound of
$$ \frac{\log \ell}{\gap} = (2/\pi^2+o(1)) \ell^4 \log \ell\ .$$

\section{Heuristic arguments for the true constants}
\label{sec:heuristic}

Up until now we have given upper bounds and lower bounds for various
mixing times and coupling times, and these bounds have typically
differed by small constant factors.  In this section we give heuristic
arguments and summarize experimental results for determining the true
asymptotic constant factors that were given in \tbref{times}.  Readers
concerned primarily with rigorous arguments will find in this section
a few theorems and many open problems.

\subsection{A million shuffles or seven}
\label{sec:million}

It is well-known that for any Markov chain, when one considers the
distance from stationarity of the distribution at time $t$, the
variation distance decays as $d(t) = (1+o(1)) A_d |\lambda|^t$, and
similarly the separation distance decays as $s(t) = (1+o(1)) A_s
|\lambda|^t$, where $\lambda$ is the second largest eigenvalue (in
absolute value).  For the Markov chains considered here, we have
rigorous exact values for $\lambda$.  To paraphrase
\cite{diaconis:cutoff}, the goal of finding mixing times is not to
determine precisely how far from stationarity a deck of cards is after
a million shuffles, but to determine if seven shuffles are enough.
For many Markov chains, the variation distance from uniformity stays
close to $1$ for a time, and then rapidly becomes small and decays
exponentially fast (see \citet{diaconis:monograph}).  The
seven-shuffle question, which is more relevant to practical
applications, asks where this cutoff occurs.  The
million-shuffle question has the virtue of typically being easier to answer,
and it appears to be relevant to the seven-shuffle question.
\cite{diaconis:cutoff} himself explains that the long-term behavior of
the Markov chain can be used as a heuristic for predicting which
Markov chains will exhibit the ``cutoff phenomenon'' in the time it
takes to randomize.  Specifically for reversible Markov chains he uses
the $L_2$ bound
$$ 4 \| P_x^t-\pi\|^2 \leq \sum_{i=1}^{N-1} v_i(x)^2 \lambda_i^{2t} $$
where the $v_i$'s are an orthogonal eigenbasis,
and uses the lead term for large $t$ 
$$\sum_{i:|\lambda_i|=\lambda} v_i(x)^2 \lambda^{2t} = \left[A_{L_2} \lambda^t\right]^2$$
to make the prediction: if $A_{L_2}$ is large then the Markov chain
probability exhibits a sharp transition.  This extends an earlier
heuristic that \citet*[sect.~7]{aldous-diaconis:walk-times} gave for
predicting the order of magnitude of the mixing time cutoff for random
walk on groups.

In this section we will work
more directly with the separation and variation distances rather than
use the $L_2$ norm, and hypothesize that
$$d(t) \approx \min\big(1,A_d |\lambda|^t\big) \ \ \ \ \text{and}\ \ \ \ 
  s(t) \approx \min\big(1,A_s |\lambda|^t\big)\ .
$$
\textit{A priori\/} it is not clear why the variation distance should
be well approximated by $A_d |\lambda|^t$ whenever this approximation
is not obviously bad, i.e.\ when it is not larger than $1$.  Indeed,
there are examples where this type of approximation fails: for walks
on random Cayley graphs on $\Z_2^d$, the variation distance and the
$L_2$ norm bound have essentially the same lead term behavior, but
exhibit sharp transitions at different times \citep{wilson:z2d}.
Nonetheless the above approximation is valid for many Markov chains
(see e.g.\ \citet*{diaconis:monograph}),
and numerical computations described in \sref{numerical} give every
indication that this approximation holds for the Markov chains that
we are interested in.  We therefore take a moment to formalize this
observation:
\begin{definition}
A family of Markov chains exhibits a \emph{clean (variation) cutoff}
if for every $\varepsilon$ there is a $K$ so that for any Markov chain
in the family and for any time $t$, whenever $|\log(A_d \lambda^t)|>K$,
$$\big|\log d(t) - \min(0,\log(A_d \lambda^t))\big| < \varepsilon\ .$$
\end{definition}
A \textit{clean separation cutoff\/} is defined similarly.  If a family
of Markov chains exhibits a clean cutoff and $A\rightarrow\infty$ then
it will necessarily exhibit a sharp cutoff (mixing time threshold)
at $\log A / \log(1/|\lambda|)$.
Since we already have the second largest eigenvalue for several classes
of Markov chains, our goal in this section
is to compute $A_s$ and $A_d$, and report on experiments that suggest
rather strongly that that the Markov chains we are considering do in
fact exhibit clean cutoffs.

\subsection{Preliminaries}

To obtain our conjectured values for the true constant factors in the
mixing times of the adjacent transposition Markov chain on
permutations and on lattice paths and the Luby-Randall-Sinclair chain
on lozenge tilings of a hexagon, we will compute $A_s$ and approximate
$A_d$ for these Markov chains.  Before working on these specific
chains, for the reader's convenience we start with some basic preliminaries
that are common to all these examples.

For the Markov chains considered here, we have a ``potential
function'' $\Phi$ (defined in \eqref{eqn:disp}, \eqref{eqn:phi-lozenge},
and \eqref{eqn:phi-sn}) such that $E[\Phi(X_{t+1})|X_t] = \lambda
\Phi(X_t)$.  If we view $\Phi$ as the vector which is $\Phi(s)$ in its
$s$th coordinate (where $s$ is a state of the chain), then $\Phi$ is
an eigenvector of the state transition matrix with eigenvalue
$\lambda$.  Since the Markov chains are monotone, and $\Phi$ is
monotone increasing with respect to this partial order, we know that
$\lambda$ is the second largest eigenvalue in absolute value.  We do
not \textit{a priori\/} know the multiplicity of the eigenvalue
$\lambda$.

Since the Markov chains considered here are reversible, their state
transition matrices are diagonalizable, and there is an orthogonal
eigenbasis $v_i$ with eigenvalues $\lambda_i$.  If the Markov chain is
started in state $s$, the distribution at time $t$ is $\sum_i \alpha_i
\lambda_i^t v_i$ where $\alpha_i = v_i(s)/(v_i\cdot v_i)$.

To determine the quantity $d(t)$ we need the worst starting state, for
$\bar{d}(t)$ we need the worst pair of starting states, and for the
separation distance $s(t)$ we need the worst start and destination
states.  It seems intuitive that for all three measures the worst
states must be the top state $\1$ and bottom state $\0$, though there
are some monotone Markov chains for which this intuition is wrong
\citep{haggstrom:personal}.
Let us consider first the Markov chain on the symmetric group.  Since
it is vertex-transitive, all
starting states are equivalent, so $\0$ and $\1$ are worst starting
states.  Furthermore, if the chain
starts in $\0$ then the worst destination state is $\1$;
\cite{fill:interruptible} used this fact in the monotone version of
Fill's algorithm.  For the other Markov chains, we can only argue that
$\0$ and $\1$ are the right states to look at for sufficiently large time
$t$, and then only under the (apparently correct) assumption that the
second largest eigenvalue $\lambda$ has multiplicity one.  While not
completely satisfying, this will allow us to compute e.g.\ $A_s$ under
this one assumption.  Assuming that $\lambda$ has multiplicity one,
from the above formula we see that when $t$ is big enough, the worst
starting states are the states maximizing $|\Phi(x)|$, i.e.\ states
$\0$ and $\1$, the worst pair of starting states are those that
maximize $|\Phi(x)-\Phi(y)|$, i.e.\ $\0$ and $\1$, and the worst start
and destination states are those that maximize $-\Phi(x)\Phi(y)$,
i.e.\ $\0$ and $\1$.

Let us suppose that the eigenvalue $\lambda$ has multiplicity one (the multiplicity is larger for the permutation chain, but we will deal with that later).
Recall that $P^t_x$ denotes the distribution at time $t$ when the chain starts in state $x$, and that $U$ denotes the uniform (stationary) distribution.
For large times $t$,
$P^{t}_{\1}-U \doteq \frac{\Phi(\1)}{\Phi\cdot \Phi}\lambda^t \Phi$.
If we view $\Phi$ as the random variable obtained by picking a uniformly random state $X$ and returning $\Phi(X)$, then
$\Phi\cdot \Phi = \Var[\Phi] N$ where $N$ is the number of states.
The vector $\Phi$ takes on its most negative value at $\0$, so it contributes $-\Phi(\0) N$ to the separation distance.  Hence for
large times $t$ we expect the separation distance to be
  $-\frac{\Phi(\1)\Phi(\0)}{\Var[\Phi]} \lambda^t$, so that
\begin{equation}
\label{eq:A_s}
A_s\meq -\frac{\Phi(\1)\Phi(\0)}{\Var[\Phi]} = \frac{\Phi(\1)^2}{\Var[\Phi]}\ ,
\end{equation}
where the question mark above the equal sign reminds us that in its
derivation we used the assumption that the second largest eigenvalue has
multiplicity one.

The contribution of $\Phi$ to the variation distance is $\frac12 N
E[|\Phi|]$.  Typically it is hard to get an analytic expression for
$E[|\Phi|]$, but heuristically it is plausible that the distribution
of $\Phi$ is Gaussian, so that $E[|\Phi|] \approx \frac{2}{\sqrt{2
\pi}} \int_0^\infty x e^{-x^2/2} \,dx\,\sqrt{\Var[\Phi]} =
\sqrt{2/\pi} \sqrt{\Var[\Phi]}$.  There are interesting Markov chains,
such as random transpositions on permutations which was analyzed by
\cite{diaconis-shahshahani:random-transpositions},
for which the principal eigenvector evaluated at a random state is very
far from being approximated by a normal distribution.  Thus in principle
this approximation should be proved for the chains that we are considering,
but we will simply assert that it is
intuitively obvious that this approximate normality holds for these chains.
Assuming this approximate normality, we have
\begin{equation}
\label{eq:A_d}
A_d \meq \frac{\Phi(\1) E[|\Phi|]/2}{\Var[\Phi]} \mapprox
\frac{\Phi(\1)}{\sqrt{2\pi \Var[\Phi]}} \meq \sqrt{A_s/(2\pi)}
\end{equation}
(where we used the approximate normality assumption in the $\mapprox$
relation, and the multiplicity-one assumption in the two $\meq$ relations).
This relation between $A_d$ and $A_s$ is consistent with the
folk-wisdom that (for reversible chains) it usually takes twice as
long for the separation distance to become small as it does for
variation distance.

\remark A notable non-reversible chain where this relation fails is
the riffle-shuffle Markov chain \citep{bayer-diaconis:dovetail}.  What
failed is the relation $\alpha_i = v_i(s)/(v_i\cdot v_i)$, which
assumes reversibility.  With the correct $\alpha_i$'s the above
heuristic reasoning gives the right thresholds for the riffle-shuffle
chain as well.

\subsection{Lozenge tilings}

For the Luby-Randall-Sinclair Markov chain on lozenge tilings of the
order $\ell$ hexagon, we used $E[(\Delta\Phi)^2]\leq O(\ell)$ to get a
bound on the variance of the height function.  While there do exist
atypical configurations for which $E[(\Delta\Phi)^2]$ is this large,
it seems that more often $E[(\Delta\Phi)^2]$ is closer to $O(1)$.  If
we substitute this bound on $E[(\Delta\Phi)^2]$ into \lref{anticonverge}, we
would get that the variance in the potential function is $O(\ell^4)$,
with typical deviations from the mean being $O(\ell^2)$.  These
heuristic bounds are in fact correct in stationarity.  It is known
that the variance in the total height is $\ell^4/4$, and more
generally $abc(a+b+c)/12$ for the hexagon with side lengths
$a,b,c,a,b,c$ \cite{blum:var}.  And as we shall see from Eq.~\ref{eq:var}
below, in stationarity the variance in $\Phi$ is $abc(a+b+c)/(2\pi^2+o(1))$.
Thus for the order $\ell$ hexagon in a stationarity, the potential
function $\Phi$ is typically within $O(\ell^2)$ of its expected value
$0$.  If (as seems likely) for each time $t$, $\Phi(X_t)$ is also
typically within $O(\ell^2)$ of its expected value, then we would get
a mixing time lower bound of $\log(\ell^3/\ell^2)/ \log(1/|\lambda|) =
(16/\pi^2+o(1)) \ell^4\log\ell$.
Intuitively the lozenge tiling is random when the distribution of its
average height is close to its stationary distribution.

\begin{theorem}
\label{thm:hex}
Consider the Luby-Randall-Sinclair lozenge tiling Markov chain on the
hexagon with side lengths $a,b,c,a,b,c$, with the tower moves parallel to
the ``$c$'' sides.
Assuming the second largest eigenvalue has multiplicity one,
$$A_s = \frac{\Phi(\1)^2}{\Var[\Phi]} = \frac{c[(a+b)^2-1]}{ab(a+b+c)}\frac{\sin^2(\pi a/(a+b))}{1-\cos(\pi/(a+b))}\ .$$
\end{theorem}
% As := (a,b,c) -> c*sin(Pi/(a+b)*a)^2/(1-cos(Pi/(a+b)))/a/b/(a+b+c)*((a+b)^2-1);
% lambda := (a,b,c) -> 1-gap(a,b,c);
% gap := (a,b,c) -> (1-cos(Pi/(a+b)))/c/(a+b-1);
% Var := (a,b,c) -> 1/4/(1-cos(Pi/(a+b)))*a*b*c*(a+b+c)/((a+b)^2-1);
\begin{proof}
The first relation is just Eq.~\eqref{eq:A_s}, which is where we use the
multiplicity-one assumption.
To compute the variance of $\Phi$ in
stationarity, we let $S_i$ denote the total height in column $i$
($-a<i<b$), and write
\begin{align*}
\Phi &= \sum_{i=-a}^b (S_i-E[S_i]) \sin\frac{\pi(i+a)}{a+b} \\
\Var[\Phi] &= \sum_{-a<i,j<b} \Cov(S_i,S_j) \sin\frac{\pi(i+a)}{a+b} \sin\frac{\pi(j+a)}{a+b}
\end{align*}
and use a formula \citep{wilson:cov} for the covariances of
the $S_i$'s
$$ \Cov(S_i,S_j) = (a+i) (b-j) \times \frac{abc(a+b+c)}{(a+b)^2((a+b)^2-1)}\ \ \ \ (i\leq j)$$
to find (with the help of Maple) that
\begin{equation}
\label{eq:var}
\Var[\Phi] = \frac{1/4}{1-\cos(\pi/(a+b))}\frac{abc(a+b+c)}{(a+b)^2-1} = (1+o(1)) \frac{abc(a+b+c)}{2\pi^2}\ .
\end{equation}
%e1:=simplify(expand(sum(sum(i*(w-j)*sin(Pi*i/w)*sin(Pi*j/w),j=i..w),i=0..w)));
%e2:=simplify(expand(sum(sum(i*(w-j)*sin(Pi*i/w)*sin(Pi*j/w),j=i+1..w),i=0..w)));
% simplify(e1+e2);
To evaluate $\Phi(\1)$ we write
$$\Phi(\1) = \sum_{i=-a}^0 \frac{(a+i)bc}{a+b}\sin(\pi(a+i)/(a+b)) +
            \sum_{i=1}^b  \frac{(b-i)ac}{a+b}\sin(\pi(a+i)/(a+b)) $$
which with the help of Maple simplifies to
\begin{equation*}
\Phi(\1) = \frac{c}{2} \frac{\sin(\pi a/(a+b))}{1-\cos(\pi/(a+b))}\ . \qedhere
\end{equation*}
% factor(simplify(expand(sum((a+i)*b*c/(a+b)*sin(Pi*(a+i)/(a+b)),i=-a..-1) + sum((b-i)*a*c/(a+b)*sin(Pi*(a+i)/(a+b)),i=0..b))));
% f := 1/2 * c * sin(Pi*a/(a+b)) / (1-cos(Pi/(a+b)));
\end{proof}

When $a=b=c=\ell$, \tref{hex} gives $A_s \mapprox 32/(3\pi^2) \ell^2$ and $A_d
\mapprox \sqrt{16/(3\pi^3)} \ell$.  Since $\gap \doteq
\pi^2/16 /\ell^4$, we estimate the separation threshold to be
$(16/\pi^2) \ell^4 \log(32/(3\pi^2) \ell^2) \doteq (32/\pi^2) \ell^4 \log
\ell$ and the variation threshold to be $(16/\pi^2) \ell^4 \log \ell$,
which matches the intuitive lower bound given above.

As a check of the $E[|\Phi|]\approx \sqrt{2/\pi} \sqrt{\Var[\Phi]}$ approximation, for the $3\times 3\times 3$ cube $E[|\Phi|]\doteq 2.872$ while $\sqrt{2/\pi} \sqrt{\Var[\Phi]}\doteq 2.892$, an error of about 1\%.  Even for the $2\times2\times2$ cube, the error between $E[|\Phi|]\doteq 1.307$ and $\sqrt{2/\pi} \sqrt{\Var[\Phi]}\doteq 1.319$ is less than 1\%.
% bpp's 2 2 2 abs-phi
% bpp's 2 2 2 sqr-phi
% for 3x3x3 got l1=29/35*3^(1/2)+51/35
%               l2=243/70*3^(1/2)+243/35
\subsection{Lattice paths}

We have effectively already computed $A_s$ and approximated $A_d$ for
lattice paths in the $a\times b$ box --- just set $c=1$ in the above
formulas for the $a\times b\times c$ lozenge tiling region.

For lattice paths we can give some additional intuition.
Consider the lattice path Markov chain on a $n/2 \times
n/2$ box.  In stationarity the height fluctation at a given site near
the center of the path will be $\Theta(\sqrt{n})$, and the
fluctuations in the potential function will be $\Theta(n^{3/2})$.
Initially the potential function is $\Theta(n^2)$, and at any given
time the fluctations in the potential function are about $O(n^{3/2})$
about its expected value.  We used these facts to obtain the lower
bound on the variation mixing time of $\log(n^2/n^{3/2}) / \log(1/|\lambda|)$.
Intuitively the path is about random when its average height is close
to its stationary distribution, which would imply that the above lower bound is tight.

\subsection{Permutations}

\begin{theorem}
\label{thm:snas}
For the random adjacent transposition Markov chain on permutations of
order $n$, assuming the second largest eigenvalue has multiplicity
$n-1$, $A_s=n-1$.
\end{theorem}

\remark It is not clear to what extent it is a coincidence that $A_s$ is the
multiplicity of $\lambda$.  This relation does not hold for the tiling
or lattice path Markov chains, but there the state spaces do not have
a group structure.  For $Z_2^n$, $A_s=n$ is the multiplicity of the
second largest eigenvalue.  But for the cycle $Z_n$, the second
largest eigenvalue has multiplicity two for $n\geq 3$, while $A_s = 2$
only for even $n>3$ and $A_s = 2\cos(\pi/n)$ for odd $n\geq 3$.

\begin{proof}[Proof of \tref{snas}]
For random adjacent transpositions on permutations, for any given card $i$
one can define an eigenvector with eigenvalue $\lambda$ based on the
location $\sigma^{-1}(i)$ of that card: $$ f_i(\sigma) = \cos [\pi
(\sigma^{-1}(i)-1/2)/n]\ .$$  There is one linear dependency amongst
these eigenvectors ($\sum_i f_i = 0$), and it appears that there are
no other eigenvectors with eigenvalue $\lambda$.
Note that
$$ f_i \cdot f_i = (n-1)! \sum_{j=1}^n \cos^2(\pi(j-1/2)/n) = (n-1)! \sum_{j=1}^n \frac{1+ \cos(2\pi(j-1/2)/n)}{2} = n!/2\ .$$
By symmetry considerations $f_i\cdot f_j = f_i\cdot f_k$ when $k\neq i\neq j$.
Since $$ 0 = f_i\cdot\left(\sum_j f_j\right) = \frac{n!}{2} + (n-1)f_1\cdot f_2$$
we have that $$ f_i\cdot f_j = -\frac{n!}{2(n-1)}$$
when $i\neq j$.

As before, to determine $A_s$ we compute the coefficients of the
$f_i$'s in the eigenbasis decomposition of $1_{\1}$.  Since there
is a linear relation amongst the $f_i$'s there will be a one-parameter
family of valid sets of coefficients -- we just need one such valid
set of coefficients.
We could be methodical and use the Gram-Schmidt procedure to extract
$n-1$ orthogonal vectors from the $f_i$'s and then use these to get
a valid set of coefficients, but the guess-and-verify method is less messy.
Consider the function
\begin{align}
\Phi &= \sum_{i=1}^n f_i(\1) f_i\ .\label{eqn:phi-sn}
\intertext{We have}
f_i \cdot \Phi &= \frac{n!}{2} \left(f_i(\1) - \frac{1}{n-1} \sum_{j\neq i} f_j(\1) \right) \notag\\
&= \frac{n!}{2} \frac{n}{n-1} f_i(\1)\ .\notag
\intertext{By comparison}
 f_i \cdot 1_{\1} &= f_i(\1)\ .\notag
\end{align}
Since the dot products with the $f_i$'s are the same (up to the constant factor $n!n/(2(n-1))$), by linear algebra we conclude that $2(n-1)/(n! n)\Phi$ has the desired coefficients.  Next we evaluate this eigenfunction at $\0$ and multiply by $-N=-n!$ to obtain $A_s$:
\begin{align*}
A_s &= -n! \frac{2 (n-1)}{n! n} \sum_{i=1}^n f_i(\1) f_i(\0) \\
 &= -2\frac{n-1}{n} \sum_{i=1}^n -\cos^2(\pi(i-1/2)/n) \\
 &= n-1\ .\qedhere
\end{align*}
\end{proof}

\subsection{Shape of the thresholds}

In Figures \ref{fig:bpp}, \ref{fig:path2}, and \ref{fig:sn}, where we present numerical data for separation and variations distances, we also plot some hypothetical asymptotic curves for the separation and variation distances, in particular
\begin{align*}
s(t) &\doteq 1-\exp({-A_s \lambda^t}) \\
d(t) &\doteq \erf\left(\frac{\sqrt{\pi}}{2} A_d \lambda^t \right) \\
\bar{d}(t) &\doteq \erf(\sqrt{\pi} A_d \lambda^t)\ ,
\end{align*}
where $\erf(x)=\int_{-x}^x e^{-t^2}\, dt\, /\sqrt{\pi}$ is the error function.
For random walk on $\Z_2^d$, the variation distance $d(t)$ was shown
to take the above form by \citet*{diaconis-graham-morrison:cutoff}.  The
intuition for why it should also hold for the Markov chains that we
are interested in is essentially the same as for their proof for
$\Z_2^d$.  When $X$ is a random state drawn from the uniform
distribution, $\Phi(X)$ is well approximated by a Gaussian.  (For
$\Z_2^d$, $\Phi$ is the number of ones.)  When $t$ is near the mixing
time threshold, $\Phi(X_t)$ should be well approximated by a Gaussian
with the same variance as $\Phi(X)$ but with mean $\Phi(X_0)
\lambda^t$.  The intuition, which was made rigorous for $\Z_2^d$ but
appears difficult to prove for the other chains, is that $\Phi$ is the
best test statistic for distinguishing $X$ from $X_t$.  The amount by
which these two Gaussians fail to overlap gives the asymptotic curve
for $d(t)$.  The curve for $\bar{d}(t)$ also follows from this
heuristic reasoning.

There does not seem to be any similar intuition for why the asymptotic
curve for $s(t)$ should be what is given above, other than that it
holds for $\Z_2^d$ and other high-dimensional product graphs, and it
appears to be a good fit at least for the tiling Markov chains.

\remark There exist Markov chains for which the asymptotic curves for
$d(t)$ are not given by the above formula.  For example,
\citet*{diaconis-fill-pitman:top-random} analyze the top-to-random
shuffle and determine the curve for $d(t)$ to be given by
an explicit piecewise-analytic formula, which in particular is
different from the above formula.  Nonetheless, the chains we are
interested in seem to have more in common with random walk on $\Z_2^d$
than they do with the top-to-random shuffle, and we are confident that
the above formulas for $d(t)$ and $\bar{d}(t)$ are the correct
asymptotic shape of the transition.

\subsection{Numerical experiments}
\label{sec:numerical}

\begin{figure}[htbp]
\centerline{\psfig{figure=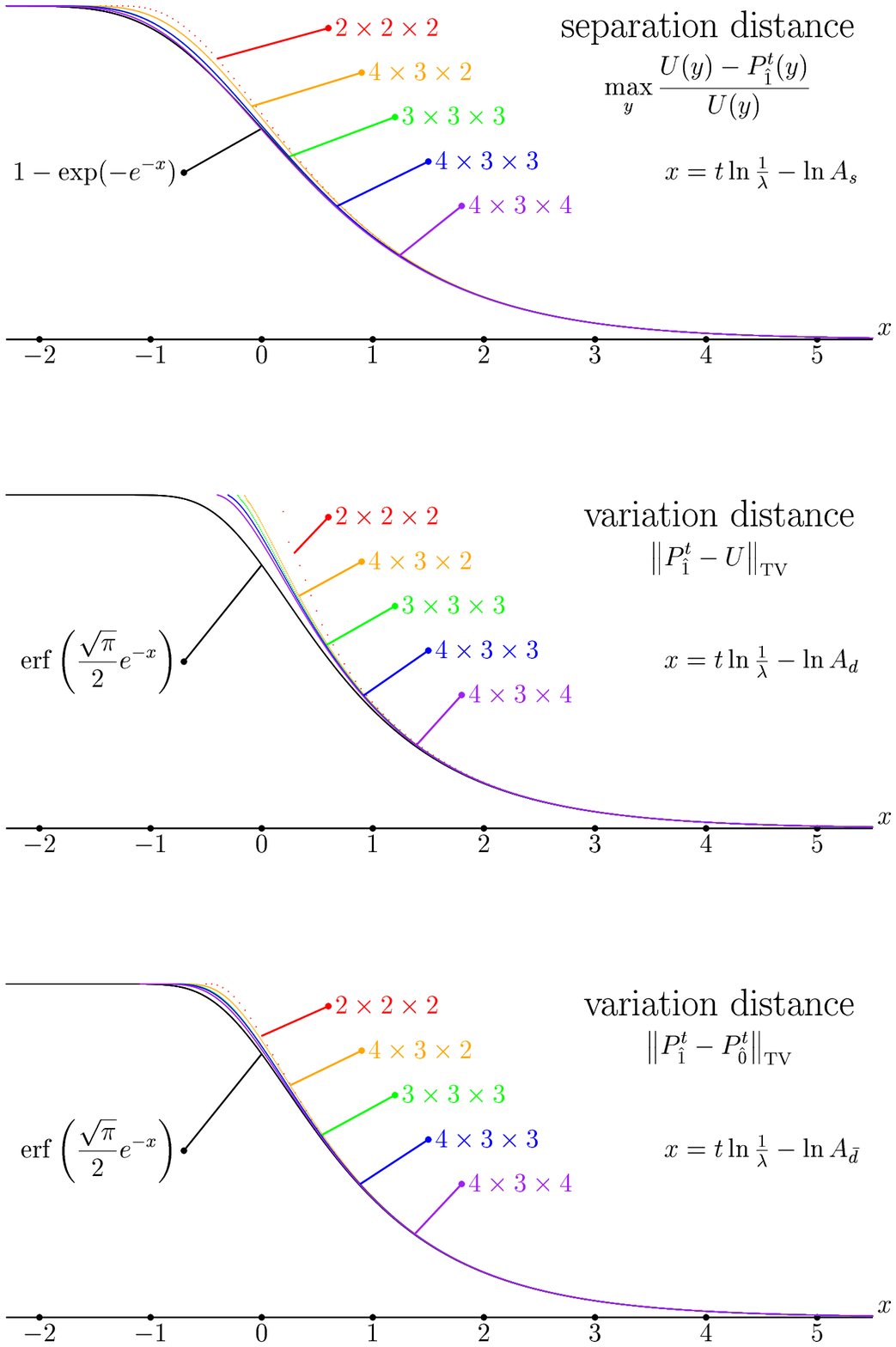}}
\caption{Data for the Luby-Randall-Sinclair lozenge tiling chain on the hexagon with side lengths $a$, $b$, $c$, $a$, $b$, $c$.}
\label{fig:bpp}
\end{figure}

\begin{figure}[htbp]
\centerline{\psfig{figure=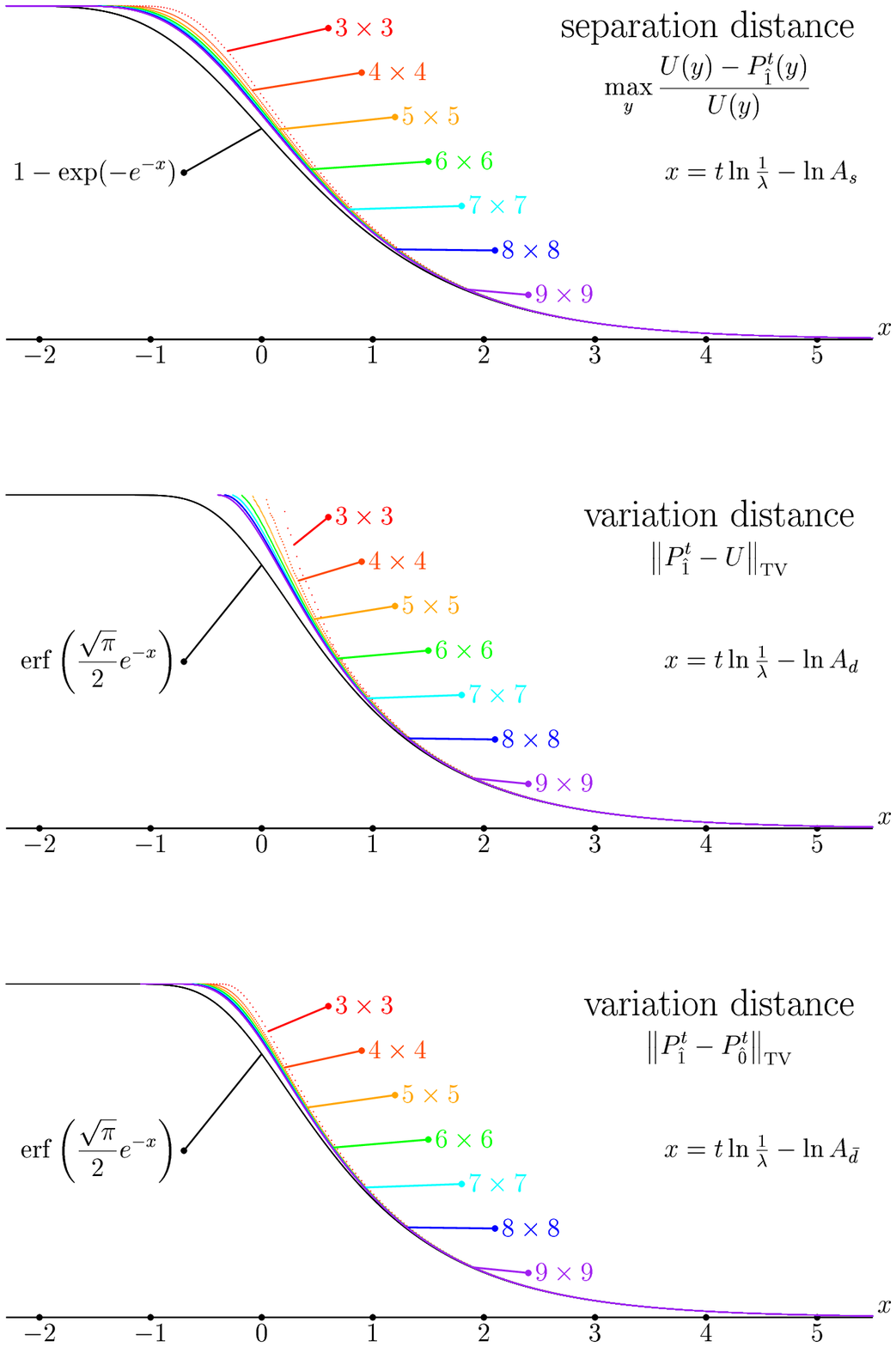}}
\caption{Data for random adjacent transpositions on $a\times b$ lattice paths.}
\label{fig:path2}
\end{figure}

\begin{figure}[htbp]
\centerline{\psfig{figure=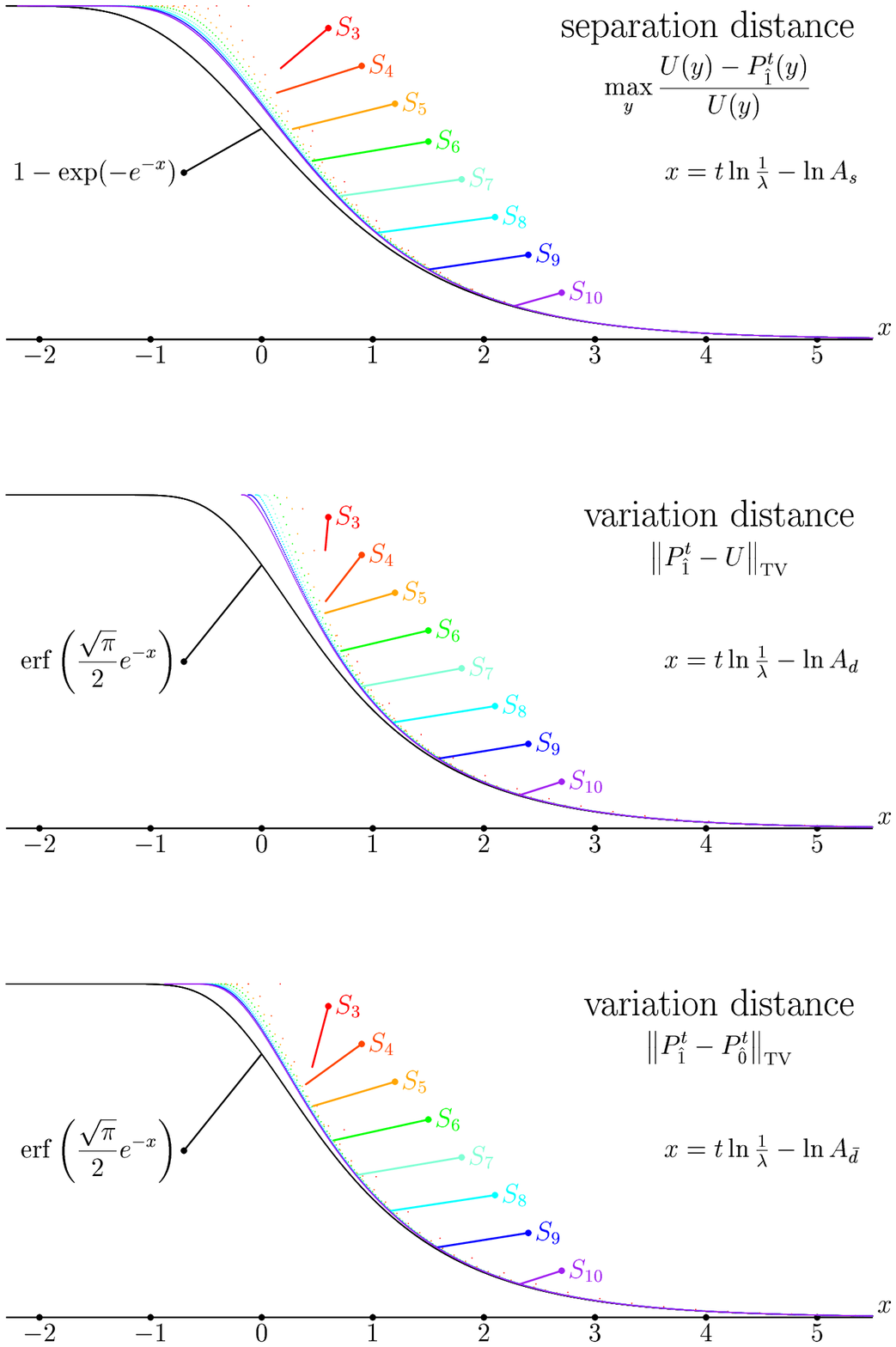}}
\caption{Data for random adjacent transpositions on the symmetric group $S_n$.}
\label{fig:sn}
\end{figure}

Figures \ref{fig:bpp}, \ref{fig:path2}, and \ref{fig:sn} show
numerical data for the convergence
rates of the three classes of Markov chains considered here.  The data
was obtained by explicit mulitplications of the state transition
matrix, and assumes that $\0$ and $\1$ are the worst starting states
for $d(t)$, worst pair of starting states for $\bar{d}(t)$, and worst
start and destination state for $s(t)$.  For each system size the
convergence data was scaled and shifted to make the transitions line
up.  The amount by which to scale and shift was computed \textit{a
priori\/} using the preceding formulas for $\lambda$, $A_s$, and
$A_d$; in particular the values for $\lambda$ and $A_s$ are exact and
the value for $A_d$ was approximated using \eqref{eq:A_d}.

The curves line up quite nicely when $d(t)$, $\bar{d}(t)$, or $s(t)$
are small, and they line up progressively better with increasing system
size.  This fact, and additional graphs not shown here, indicate that
the Markov chains have a clean cutoff, and that hence
the cutoff phenomenon does indeed
occur where we expect it.  The data looks less good when $d(t)$,
$\bar{d}(t)$, or $s(t)$ are large; this distortion is due in part to
the fact that for finite system sizes there is a small finite $x_{\min}$,
while in the idealized limit $x_{\min}=-\infty$.  For $s(t)$, $x_{\min}$
is about twice as large as for $d(t)$, so this distortion is much less
pronounced for $s(t)$.  For $\bar{d}(t)$, $x_{\min}$ is $\log 2$ more
negative than for $d(t)$, which is large enough to make the curves
line up noticeably better for $\bar{d}(t)$ than for $d(t)$.

\subsection{Monte Carlo experiments}

To estimate the coupling time for the three classes of Markov chains,
one can actually run the Markov chain until the upper and lower
configurations coalesce, repeat many times, and compare the results
for different system sizes.  The obvious advantage of Monte Carlo over
numerical experiments is that one can do much larger system sizes.
But even with the large system sizes, when comparing different system
sizes it is still much better to rescale time by
 $$\log \frac{1}{\lambda} = \log \frac{1}{1-(1-\cos(\pi/n))/(n-1)}$$
rather than its asymptotic value $\frac{\pi^2}{2} / n^3$.

After rescaling time in this manner it becomes quite clear that the
coupling time for $n/2 \times n/2$ lattice paths is about $\log n /
\log 1/\lambda \doteq 2/\pi^2 n^3 \log n$, and the coupling time for
permutations is $\log n^2 / \log 1/\lambda \doteq 4/\pi^2 n^3 \log n$.
Of course this coupling time estimate for lattice paths is actually
rigorous thanks to Theorems~\ref{thm:path-perm-upper} and \ref{thm:path-lower-c}.
Estimating the correct constant for the tiling Markov chain is however
much more challenging, and we have not yet succeeded in doing this.

Surprisingly one can use Fill's algorithm to do similar Monte Carlo
experiments to measure the separation distance
\citep[sect.~9]{fill:interruptible}; we did not do this since we
already had the numerical data.  We are unaware of any similar Monte
Carlo method for measuring the variation distance.

\old{
Next consider the gap between an upper and lower height function
(conditioned to be different) after a large number of steps.  It seems
plausible that this conditional gap has size $\Theta(\ell^2)$, which
would imply that the coupling time threshold is again
$\log(\ell^3/\ell^2) / \log(1/|\lambda|)$.
}

\section{Concluding remarks}\label{sec:conclude}

\nib
Adding weights to distance functions can be useful when proving mixing
time upper bounds.  While the optimal weighting scheme will be related
to an eigenvector of the state transition matrix, there is no need to
diagonalize the matrix, nor is it even necessary to exhibit a single
eigenvector to produce an effective weighting scheme that yields good
upper bounds.  For example, one could have used parabolic weights
rather than the sinusoidal weights that we did, and still derived
mixing time upper bounds that are only a constant factor worse than
the ones we derived.

\nib
There are a variety of Markov chains for which a mixing time cutoff
phenomenon has been proved.  In future work on sharp mixing time
thresholds, it would be worthwhile to determine whether or not the
Markov chains exhibit a clean cutoff as defined in \sref{million}.

\section*{Acknowledgements}

The author thanks Persi Diaconis and David Aldous for useful
discussions, and Jim Fill, Chris Henley, and a referee for their
detailed comments on an earlier draft.

\nocite{ng:thesis}

\bibliography{nlm}
\bibliographystyle{abbrvnat}

\end{document}